\documentclass[english]{article}
\usepackage[T1]{fontenc}
\usepackage[utf8]{inputenc}
\usepackage{enumitem}
\usepackage{amsmath}
\usepackage{amsthm}
\usepackage{amssymb}
\usepackage{esint}
\usepackage{orcidlink} 
\usepackage{tocloft}
\makeatletter
\theoremstyle{plain}
\newtheorem{thm}{\protect\theoremname}[section]
\theoremstyle{definition}
\newtheorem{example}[thm]{\protect\examplename}
\theoremstyle{plain}
\newtheorem{lem}[thm]{\protect\lemmaname}
\theoremstyle{remark}
\newtheorem{rem}[thm]{\protect\remarkname}
\theoremstyle{definition}
\newtheorem{defn}[thm]{\protect\definitionname}
\theoremstyle{plain}
\newtheorem{cor}[thm]{\protect\corollaryname}

\usepackage{hyperref}
\hypersetup{
	colorlinks=true,
	linkcolor=black,    
	citecolor=black,    
	urlcolor=blue,      
}
\usepackage{nameref}
\usepackage{mathtools}
\usepackage{tikz} 
\usetikzlibrary{positioning} 
\usepackage{multicol}
\mathtoolsset{showonlyrefs} 
\makeatletter
\let\orgdescriptionlabel\descriptionlabel
\renewcommand*{\descriptionlabel}[1]{%
  \let\orglabel\label
  \let\label\@gobble
  \phantomsection
  \edef\@currentlabel{#1}%
  \let\label\orglabel
  \orgdescriptionlabel{#1}%
}
\makeatother

\usepackage{subcaption}

\makeatother

\usepackage{babel}
\providecommand{\corollaryname}{Corollary}
\providecommand{\definitionname}{Definition}
\providecommand{\examplename}{Example}
\providecommand{\lemmaname}{Lemma}
\providecommand{\remarkname}{Remark}
\providecommand{\theoremname}{Theorem}

\begin{document}

\title{An Introduction to the Obstacle Problem}

\author{Andreas Minne~\orcidlink{0000-0002-9608-3984}
\thanks{%
CMUC, Department of Mathematics, University of Coimbra, 3000-143 Coimbra, Portugal. \newline
Email: \href{mailto:andreas.minne@mat.uc.pt}{\nolinkurl{andreas.minne@mat.uc.pt}}}
\thanks{To appear in \emph{Modern Methods in the Analysis of Free Boundary Problems}, edited by E.~A.~Pimentel and J.~M.~Urbano, Coimbra Mathematical Texts, Springer.}}

\maketitle
\global\long\def\dist{\operatorname{dist}}
\global\long\def\trace{\operatorname{tr}}
\global\long\def\interior{\operatorname{int}}
\global\long\def\df{\coloneqq}
\global\long\def\fd{\eqqcolon}

\begin{abstract}
We give a self-contained account of the classical obstacle problem, written for graduate students who wish to enter the field of free boundary problems.
Starting from a one-dimensional model and an explicit radially symmetric example, we set up the variational problem, show that minimizers exist and are unique, and prove interior $W^{2,p}$-regularity as well as the optimal $C^{1,1}$-regularity.
The second half of the notes is devoted to the free boundary: we prove that it is porous, introduce blow-up limits, and use the Alt--Caffarelli--Friedman monotonicity formula to show that every free boundary point is either regular or singular.

These notes grew out of a mini-course at the Rio Summer School on PDEs in 2020, and they are meant to be read alongside, not instead of, the monograph by Petrosyan, Shahgholian, and Uraltseva.


\end{abstract}
\tableofcontents

\section{Introduction}

Free boundary problems (FBPs) form one of the most active areas in the analysis of partial differential equations (PDEs).
They draw on classical PDE theory, functional analysis, geometric measure theory, and topology, and they arise in applications ranging from physics and biology to finance and engineering.

Rather than surveying the whole field, these notes treat a single problem in detail: the \emph{obstacle problem}, arguably the best-known free boundary problem of all.
One looks for a function that minimizes the Dirichlet energy among all functions lying above a given obstacle; the free boundary is the interface where the minimizer leaves the obstacle.
Our aim is that the reader comes away with ideas and techniques that carry over to other FBPs and related questions.

The material is based on a mini-course given at the Rio Summer School on PDEs in Rio de Janeiro in 2020, organized jointly by PUC-Rio, UFF, and UFRJ. Much of it follows the monograph by Petrosyan, Shahgholian, and Uraltseva \cite{MR2962060}, which we recommend keeping at hand while reading.

The reader is assumed to know the standard tools of elliptic PDEs and functional analysis, such as Sobolev spaces and the maximum principle.
Good references for these are the books by Gilbarg and Trudinger \cite{gilbarg:01} and by Evans \cite{evans2010partial}.

\subsection*{What Makes a Boundary Free?}

To see what distinguishes a free boundary problem, think of heat flowing through a wooden table.
Write $u$ for the temperature, $\Omega$ for the table, and $\partial \Omega$ for its boundary.
Apart from a negligible thermal expansion, the table keeps its shape, so $\Omega$ is a fixed domain and $u$ is the only unknown.

Suppose instead that the table is made of ice.
As it melts, its shape changes over time, and the temperature $u$ depends on how $\partial \Omega$ moves.
At the same time, the motion of $\partial \Omega$ depends on $u$, since the parts of the table that receive more heat melt faster.

So for the wooden table only $u$ has to be found, while for the ice table \emph{both} $u$ and $\partial \Omega$ are unknown.
This is the essential difference between a PDE in a fixed domain and a free boundary problem.
The coupling between $u$ and $\partial \Omega$ makes the analysis considerably harder, and understanding the geometry and regularity of $\partial\Omega$ is one of the central challenges.

The notes are organized as follows. \S~\ref{sec:Setup-and-Notation} collects the notation. \S~\ref{sec:main-section} builds intuition from a one-dimensional model. \S~\ref{sec:general-case} treats a radially symmetric example in higher dimensions and then develops existence and optimal regularity in general. \S~\ref{sec:free-boundary} is devoted to the free boundary.

\section{Notation\label{sec:Setup-and-Notation}}

The following notation is mostly standard but is included here for completeness and easy reference.

\begin{multicols}{2}
	\begin{itemize}[itemsep=0pt,leftmargin=*]
	\item
$\mathcal{D}\subseteq\mathbb{R}^{n}$ is a bounded open set with $C^{1}$ boundary.
	\item
$B_{r}(x)$ denotes the ball in $\mathbb{R}^{n}$ of radius $r$ centered at $x$; we also write $B_{r}$ if $x$ is the origin.
	\item
$x_1$, $x_2$ etc.~are coordinates for $x\in\mathbb{R}^n$.
	\item
$x^0$, $x^1$ etc.~denote points in $\mathbb{R}^n$.
	\item
$C^{k,\alpha}(\mathcal{D})$ denotes the space of functions whose derivatives up to order $k$ are $\alpha$-H\"{o}lder continuous (or Lipschitz if $\alpha=1$) in $\mathcal{D}$.
	\item
$L^p(\mathcal{D})$ are the classical Lebesgue function spaces.
	\item
$W^{k,p}(\mathcal{D})$ and $W_0^{k,p}(\mathcal{D})$ denote the standard Sobolev spaces, the latter with zero traces.
	\item
$u$ is the solution to the obstacle problem.
	\item
$\psi$ is the obstacle function (and later set to the zero function).
	\item
$\nabla u$ and $D^2 u$ denote the gradient and Hessian of $u$ respectively.
	\item
$\Delta\df\sum_{k=1}^{n}\frac{\partial^2}{\partial x_{k}^{2}}$ is the Laplace operator.
	\item
The trace of a matrix $M$ is denoted by $\trace {M}$, hence $\Delta u = \trace(D^2 u)$.
	\item
$\Omega(u)\df\{x\in\mathcal{D}:u(x)>\psi(x)\}$ is the noncoincidence set of $u$.
Also the notation $\{u>\psi\}$ is used.
	\item
$\Lambda(u)\df\{x\in\mathcal{D}:u(x)=\psi(x)\}$ is the coincidence set of $u$.
	\item
$\Gamma(u) \df \partial \Omega(u) \cap \mathcal{D}$ is the free boundary.
	\item
$\mathcal{K}_{\psi}\df\{v\in W_{0}^{1,2}(\mathcal{D}):v\ge\psi\text{ a.e.}\}$.
	\item
$f^{+}\df\max(f,0)$, $f^{-}\df\max(-f,0)$.
	\item
For a set $E$, $\chi_{E}(x)\df\begin{cases} 1 & x\in E,\\ 0 & x\in E^c. \end{cases}$
	\item
$B\Subset\mathcal{D}$ means that $\overline{B}$ is compactly contained in $\mathcal D$.
	\item
$\lvert\mathcal{D}\rvert$ is the Lebesgue measure of $\mathcal{D}$.
	\item
$u_{k}\rightharpoonup u$ denotes weak convergence of $\{u_k\}_{k=1}^\infty$ to $u$.
	\item
$u_{k}\to u$ denotes strong convergence of $\{u_k\}_{k=1}^\infty$ to $u$.
	\item
The constant $C$ may vary from line to line in estimates, as is standard in PDE literature.
	\item
$\fint_{\mathcal{D}} f(x)\,dx \df \frac{1}{\lvert\mathcal{D}\rvert}\int_{\mathcal{D}} f(x)\,dx$.
	\item
$\interior(E)$ denotes the open set consisting of interior points of the set $E$.
\end{itemize}
\end{multicols}

\section{A One-Dimensional Model: The Rubber Band\label{sec:main-section}}

\begin{figure}
\begin{minipage}[t]{0.45\textwidth}%
\centering
\begin{tikzpicture}
\node[anchor=south east,inner sep=0] (image) at (0,0) {\includegraphics[width=\textwidth]{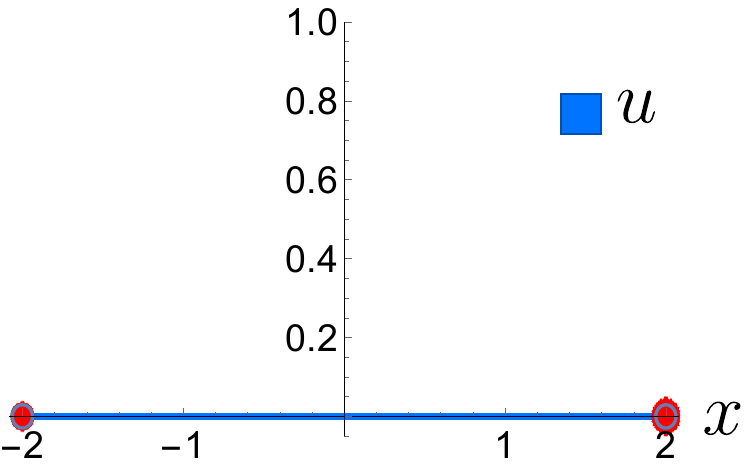}};
\node[below left] (copyright-text) {\tiny{\copyright Andreas Minne}};
\end{tikzpicture}\subcaption{Rubber band without an obstacle.}\label{fig:1Dslnwithoutobstacle}%
\end{minipage}\hfill{}%
\begin{minipage}[t]{0.45\textwidth}%
\centering
\begin{tikzpicture}
\node[anchor=south east,inner sep=0] (image) at (0,0) {\includegraphics[width=\textwidth]{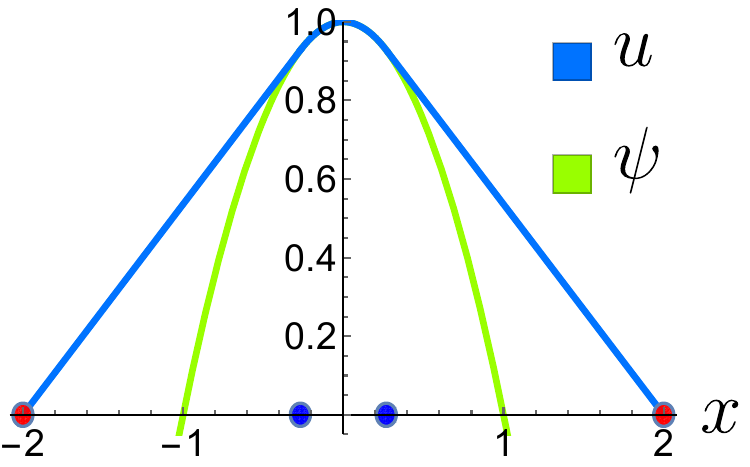}};
\node[below left] (copyright-text) {\tiny{\copyright Andreas Minne}};
\end{tikzpicture}\subcaption{Rubber band with an obstacle.}\label{fig:1Dslnwithobstacle}%
\end{minipage}

\caption{}
\end{figure}

Imagine a rubber band held under tension.
In this case, the gravitational force is negligible and the band --- which minimizes the energy --- looks like a line segment, as illustrated by the graph of $u$ in Figure~\ref{fig:1Dslnwithoutobstacle}; the coordinate system has been chosen so that the ends of the rubber band are fixed at the coordinates $(-2,0)$ and $(2,0)$, and $u(x)$ is the height of the band measured from the $x$-axis.

However, if there is an obstacle $\psi$ between the end points (which in our case looks suspiciously similar to $1-x^{2}$), the rubber band would instead look like the graph of $u$ in Figure~\ref{fig:1Dslnwithobstacle}.

At this point the attentive reader might wonder, ``So where is the free boundary?''
Good question!
The answer is given by the two points at which the rubber band goes from touching the obstacle to not touching it; these two points are highlighted as two blue dots on the $x$-axis in Figure~\ref{fig:1Dslnwithobstacle}, and in Example~\ref{exa:1D} the free boundary points will be explicitly calculated.
Note that, unlike the table made of ice discussed in the introduction, the free boundary in this problem is not free in the sense that it can move (this is a stationary problem after all) --- instead it is free in the sense that it is \emph{unknown a priori}.

The energy that the rubber band minimizes is proportional to the length functional (since gravity is negligible when the band is under tension), which we in turn approximate with a Taylor expansion to simplify the analysis, \[ L(u)\df\int_{-2}^{2}\sqrt{1+(u'(x))^{2}}\,dx\approx\int_{-2}^{2}\Big(1+\frac{(u'(x))^{2}}{2}\Big)\,dx.
\]

This Taylor expansion would normally let us turn to a linear PDE instead of the nonlinear minimal surface equation; however, we will see that in the obstacle problem the PDE stays nonlinear regardless, although the simplification above is still helpful.
Hence $u$ minimizes the Dirichlet energy $\int_{-2}^{2}(u'(x))^{2}\,dx$ in the set \[ \mathcal{K}_{\psi}\df\{v\in W_{0}^{1,2}(-2,2):v(x)\ge\psi(x)\text{ a.e.}\text{ in }(-2,2)\}, \] where $W_{0}^{1,2}(\mathcal{D})$ denotes the standard Sobolev space for which the functions are zero on the boundary $\partial\mathcal{D}$ in the sense of traces, and we recall that $\psi$ denotes the obstacle.
Since we are currently working in one dimension, $v$ is $1/2$-H\"{o}lder continuous.
We also introduce the convenient notation $\Omega(u)\df\{x\in\mathcal{D}:u(x)>\psi(x)\}$ for the \emph{noncoincidence set} and $\Lambda(u)\df\{x\in\mathcal{D}:u(x)=\psi(x)\}$ for the \emph{coincidence set}.
The \emph{free boundary} is denoted by $\Gamma(u)\df\partial\Omega(u)\cap\mathcal{D}$.\footnote{Later on, we will see that the obstacle $\psi$ can be set to 0, so $\Omega(u)=\{x\in\mathcal{D}:u(x)>0\}$, $\Lambda(u)=\{x\in\mathcal{D}:u(x)=0\}$, and $\Gamma(u)=\partial\Omega\cap\mathcal{D}$.}

Since $u$ is a minimizer and continuous, a standard perturbation argument with $u+\epsilon\phi$, where $\phi$ is supported on the noncoincidence set $\Omega(u)$, shows that $u$ is harmonic in $\Omega(u)$.
To get acquainted with the problem, let us flesh out the details in the example below, with some help from the heuristics given by Figure~\ref{fig:1Dslnwithobstacle}.
\begin{example}
\label{exa:1D} Because of symmetry of the problem, we only need to
consider nonnegative values of $x$.
Let $x_{f}$ be the ($x$-coordinate of the) free boundary point where the rubber band leaves the obstacle, indicated by a blue dot on the positive $x$-axis in Figure~\ref{fig:1Dslnwithobstacle}.
From experience, we know that the solution looks like $u_{1}$ in Figure~\ref{fig:1D-half-solution}, and not like $u_{2}$ for which the free boundary point $x_{f}$ is further to the right.
The difference is that $u_{1}$ leaves the obstacle $\psi$ tangentially while $u_{2}$ does not, i.e.,~we have the condition $u'(x_{f})=\psi'(x_{f})$ ($u_{2}$ does not even have a well-defined derivative at its free boundary point).
\begin{figure}
\centering
\begin{tikzpicture}
\node[anchor=south east,inner sep=0] (image) at (0,0) {\includegraphics[width=0.65\textwidth]{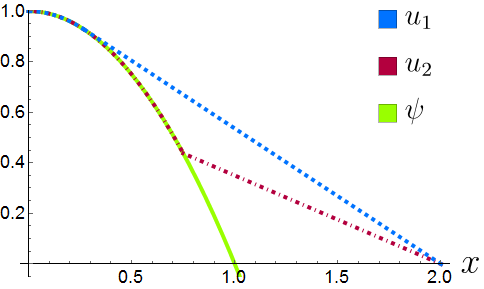}};
\node[below left] (copyright-text) {\tiny{\copyright Andreas Minne}};
\end{tikzpicture}\caption{The rubber band leaves the obstacle $\psi$ tangentially, hence $u_{1}$
(and not $u_{2}$) is the correct solution.\label{fig:1D-half-solution} }
\end{figure}

From an energetic perspective we can also justify the condition $u'(x_{f})=\psi'(x_{f})$ since $u$ can be thought of as ``wanting to leave'' the (nonharmonic) obstacle to become harmonic as soon as possible.
Yet it cannot leave earlier than tangentially since it would then cross the obstacle.

If we collect all information, and recall that $u$ is harmonic in the noncoincidence set $\Omega=(x_{f},2)$, then \[
\begin{cases}
u(x)=\psi(x) & \text{in }\Lambda=[0,x_{f}],\\
u''(x)=0 & \text{in }\Omega=(x_{f},2),\\
u(x_{f})=\psi(x_{f}),\\
u'(x_{f})=\psi'(x_{f}),\\
u(0)=1,\\
u(2)=0.
\end{cases}
\] From the second equation, $u(x)=ax+b$ in $\Omega(u)$, and $b=-2a$ from the condition $u(2)=0$.
Moreover, $u'(x_{f})=a=\psi'(x_{f})=-2x_{f}$.
If we plug this into the equation $u(x_{f})=\psi(x_{f})$ we get $-2x_{f}^{2}+4x_{f}=1-x_{f}^{2}$ from which $x_{f}=2-\sqrt{3}$ is the only viable root on the positive $x$-axis.
Hence we can summarize our result: The coincidence set $\Lambda=[0,2-\sqrt{3}]$, the free boundary $\Gamma=\{2-\sqrt{3}\}$ and the noncoincidence set $\Omega=(2-\sqrt{3},2)$.
The solution is \[ u=\begin{cases} 1-x^{2}, & x\in[0,2-\sqrt{3}],\\ (2\sqrt{3}-4)x+4(2-\sqrt{3}), & x\in(2-\sqrt{3},2].
\end{cases}
\] A corresponding calculation can of course be made for negative values of $x$.
\end{example}

There are two noteworthy remarks that we can have in mind for the future analysis of the obstacle problem.

The first is that since we have more unknowns (both $u$ and $\Gamma$) than in a traditional PDE, we need more information/conditions to uniquely solve the problem.
In this smooth 1D setting, these conditions manifest themselves in the free boundary conditions $u|_{\Gamma}=\psi|_{\Gamma}$ and $u'|_{\Gamma}=\psi'|_{\Gamma}$.

The second remark is that even though the mathematical ingredients --- such as the Laplacian and the obstacle (in this example) --- are analytic, $u$ is \emph{not} analytic.
In fact, $u$ is not even $C^{2}$ since the second derivatives are discontinuous across the free boundary (recall that $x_{f}=2-\sqrt{3}$ is the free boundary point in the example),
\begin{align*}
\lim_{x\to x_{f},x>x_{f}} & u''(x)=0,\\
\lim_{x\to x_{f},x<x_{f}} & u''(x)=-2.
\end{align*}
This is a characteristic feature of solutions to the obstacle problem.
However, the first derivatives are Lipschitz, so $u\in C^{1,1}(-2,2)$, which consequently is the best regularity one could hope for.

\section{The Obstacle Problem in Higher Dimensions}\label{sec:general-case}

Now that we have dipped our toes into the obstacle problem by studying the one-dimensional case, let us turn to the corresponding setting in $\mathbb{R}^{n}$, $n\ge2$.
Instead of a rubber band, we have a (higher-dimensional) membrane in $\mathbb{R}^{n}$ that has some fixed boundary condition and which has to stay over some given obstacle.

\subsection{A Radially Symmetric Example}
\begin{figure}
\centering
\begin{tikzpicture}
\node[anchor=south east,inner sep=0] (image) at (0,0) {\includegraphics[width=0.65\textwidth]{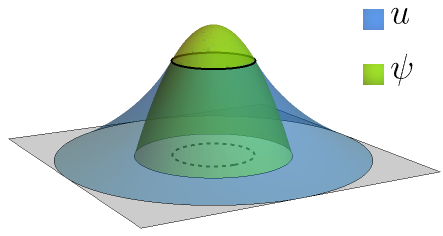}};
\node[below left] (copyright-text) {\tiny{\copyright Andreas Minne}};
\end{tikzpicture}\caption{Due to radial symmetry, the free boundary is a circle. \label{fig:2dsln}}
\end{figure}

\begin{figure}
\centering
\begin{tikzpicture}
\node[anchor=south east,inner sep=0] (image) at (0,0) {\includegraphics[width=0.65\textwidth]{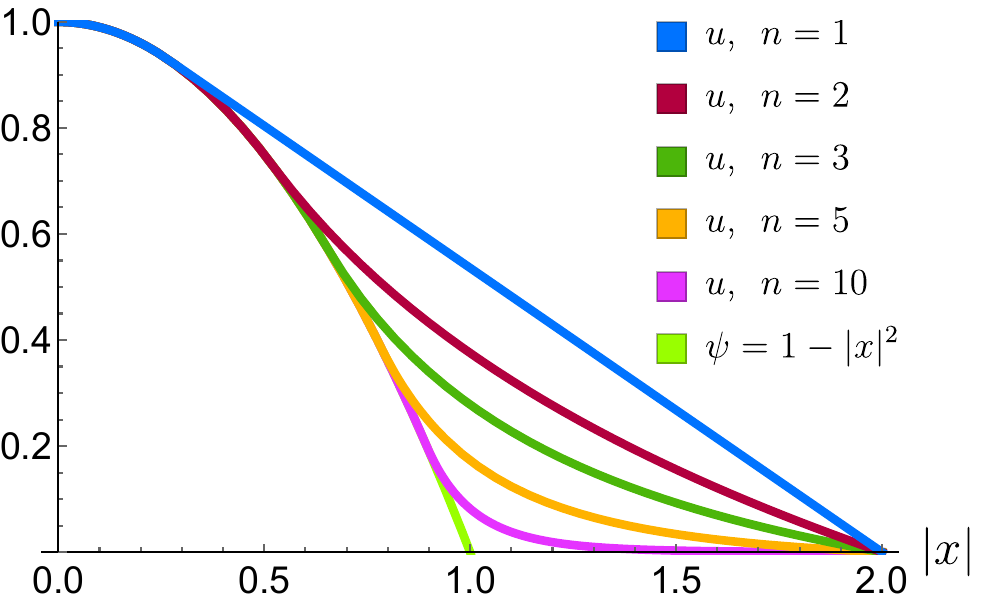}};
\node[below left] (copyright-text) {\tiny{\copyright Andreas Minne}};
\end{tikzpicture}

\caption{The free boundary radius $r_{f}$ increases with dimension. \label{fig:slns-many-dims}}
\end{figure}

To ease into this setting, let us consider a radial configuration, where the fixed domain is $\mathcal{D}=B_{2}$ with $u=0$ on the boundary $\partial B_{2}$, and the obstacle $\psi(x)=1-\lvert x\rvert ^{2}$.
The membrane looks like the graph of $u$ in Figure~\ref{fig:2dsln}, and the free boundary is naturally a circle (with a radius yet to be determined).
As in the one-dimensional case, from the Dirichlet energy functional we derive and solve the following equations: \[
\begin{cases}
u=\psi & \text{in }\Lambda,\\
\Delta u=0 & \text{in }\Omega,\\
u|_{\Gamma}=\psi|_{\Gamma}, & \nabla u|_{\Gamma}=\nabla\psi|_{\Gamma},\\
u|_{\partial B_{2}}=0.
\end{cases}
\] By symmetry, the solution $u$ is radial, i.e., $u(x) = u(r)$ with $r = \lvert x\rvert$.
Hence \[ \Delta u=\frac{1}{r^{n-1}}\frac{d}{dr}\left(r^{n-1}\frac{du}{dr}\right)=0 \] in the noncoincidence set $\Omega$.
Similar calculations as in Example~\ref{exa:1D} show that the radius $r_{f}$ defining the free boundary solves
\begin{align*}
2r_{f}^{2}\ln r_{f}-(2\ln2+1)r_{f}^{2}+1 & =0,\qquad n=2,\\
2^{-(n-3)}r_{f}^{n}-nr_{f}^{2}+n-2 & =0,\qquad n\ge3.
\end{align*}
In two dimensions ($n=2$) we get $r_{f}\approx0.52$ while $r_{f}\approx0.65$ for $n=3$.
Solutions for different dimensions are graphed as a function of the radius $r=\lvert x\rvert$ in Figure~\ref{fig:slns-many-dims}.
We observe that $r_{f}$ monotonically increases and converges to one as $n\to\infty$.

\subsection{Variational Formulation}

The obstacle problem asks for the shape of a membrane constrained to lie above a given obstacle.
By modeling the membrane's energy and neglecting higher-order effects (e.g.,~large slopes or gravity), we arrive at the Dirichlet energy as the functional to minimize.
Having seen this in the one-dimensional and radial cases, we now write the problem in its general variational form, for a domain $\mathcal{D} \subseteq \mathbb{R}^n$: \[ \min_{\mathcal{K}_{g,\psi}}\int_{\mathcal{D}}\lvert\nabla u\rvert^{2}\,dx, \] where $\mathcal{K}_{g,\psi}\df\{u\in W^{1,2}(\mathcal{D}):u-g\in W_{0}^{1,2}(\mathcal{D}),u\ge\psi\text{ a.e.~in }\mathcal{D}\}$.
The boundary condition $g$ is assumed to be in $W^{1,2}(\mathcal{D})$, the obstacle $\psi\in C^{2}(\mathcal{D})$, and $\mathcal{D}$ is a bounded $C^{1}$ domain.
We also impose the compatibility condition that the obstacle does not interfere with the fixed boundary condition, i.e.,~$g\ge\psi$ on $\partial\mathcal{D}$ with the stringent meaning that $(\psi-g)^{+}\in W_{0}^{1,2}(\mathcal{D})$ (we use the notation $f^{+}\df\max(f,0)$ and $f^{-}\df\max(-f,0)$).

\paragraph*{Reduction to the case of zero obstacle}

There is a standard procedure to recast the minimizing problem above as an equivalent problem in which the obstacle is removed in favor of a forcing term.
The details are spelled out below.

\begin{figure}
\begin{tabular}{p{5cm}c}
{\begin{equation*}
\begin{cases}
\Delta u=0 & \text{in }\Omega\df\{u>\psi\},\\ 
\phantom{\Delta}u=\psi & \text{in }\Omega^{c},\\
\nabla u=\nabla\psi & \text{on }\partial\Omega.
\end{cases}
\end{equation*}}
&
\raisebox{-80pt}{\includegraphics[width=0.45\textwidth]{Figures/FB1D.pdf}}
\end{tabular}

\medskip
Renaming $u-\psi$ as $u$ and setting $f\df-\Delta\psi$,
\medskip

\begin{tabular}{p{5cm}c}
{\begin{equation*}
\begin{cases}
\Delta u=f & \text{in }\Omega\df\{u>0\},\\ 
\phantom{\Delta}u=0 & \text{in }\Omega^{c},\\
\nabla u=0 & \text{on }\partial\Omega.
\end{cases}
\end{equation*}}
&
\raisebox{-90pt}{\begin{tikzpicture}
\node[anchor=south east,inner sep=0] (image) at (0,0) {\includegraphics[width=0.45\textwidth]{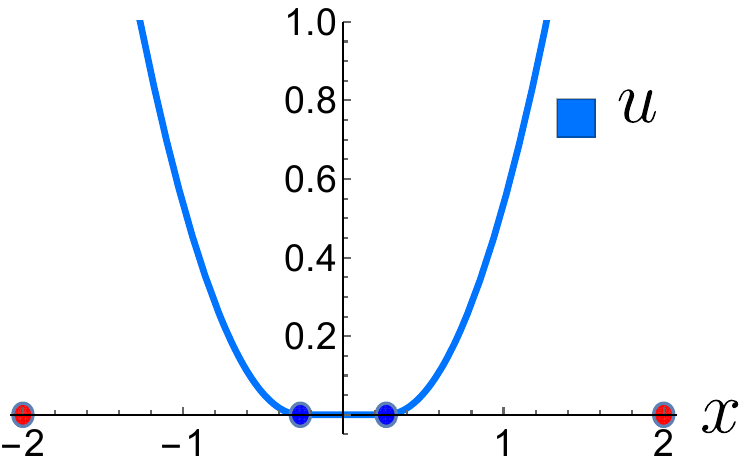}};
\node[below left] (copyright-text) {\tiny{\copyright Andreas Minne}};
\end{tikzpicture}}
\end{tabular}

\caption{Reduction to zero obstacle.\label{fig:reduction-zero-obstacle}}
\end{figure}

Let $u$ be a minimizer, $\psi$ the obstacle and define $v$ as $v\df u-\psi$.
\begin{align*}
\int_{\mathcal{D}}\lvert\nabla u\rvert^{2}\,dx & =\int_{\mathcal{D}}\lvert\nabla(v+\psi)\rvert^{2}\,dx=\int_{\mathcal{D}}\lvert\nabla v\rvert^{2}+2\nabla v\cdot\nabla\psi+\lvert\nabla\psi\rvert^{2}\,dx\\
 & =\{\text{IBP}\}=\int_{\mathcal{D}}\lvert\nabla v\rvert^{2}-2v\Delta\psi+\lvert\nabla\psi\rvert^{2}\,dx+2\int_{\partial\mathcal{D}}v\frac{\partial\psi}{\partial\nu}\,dS\\
 & =\int_{\mathcal{D}}\lvert\nabla v\rvert^{2}-2v\Delta\psi+\lvert\nabla\psi\rvert^{2}\,dx+2\int_{\partial\mathcal{D}}(g-\psi)\frac{\partial\psi}{\partial\nu}\,dS
\end{align*}
where the boundary integrands are taken in the trace sense.
Let $f\df-\Delta\psi$.
Since only the first two terms depend on $v$, it follows that $v$ minimizes the functional \[ J(w)=\int_{\mathcal{D}}\lvert\nabla w\rvert^{2}+2wf\,dx \] in $\mathcal{K}_{g_{1},0}\df\{w\in W^{1,2}(\mathcal{D}):w-g_{1}\in W_{0}^{1,2}(\mathcal{D}),w\ge0\text{ a.e.~in }\mathcal{D}\}$, where $g_{1}\df g-\psi$.
Geometrically, this substitution is demonstrated in Figure~\ref{fig:reduction-zero-obstacle}.

\subsection{Existence, Uniqueness, and \texorpdfstring{$W^{2,p}$}{W2,p}-Regularity}

Let us briefly summarize the current setting: We seek minimizers of \[ J(v)\df\int_{\mathcal{D}}\lvert\nabla v\rvert^{2}+2vf\,dx \] over the convex set \[ \mathcal{K}_{g,0}\df\{v\in W^{1,2}(\mathcal{D}):v-g\in W_{0}^{1,2}(\mathcal{D}),v\ge0\text{ a.e.~in }\mathcal{D}\}, \] where $g\in W^{1,2}(\mathcal{D})$ satisfies $g\ge0$ on $\partial\mathcal{D}$ in the sense that $g^{-}\in W_{0}^{1,2}(\mathcal{D})$, and $f\in L^{\infty}(\mathcal{D})$.\footnote{For notational simplicity, we write $g$ instead of $g_1$.}
In fact, since $f$ stems from $\Delta\psi$ for $\psi\in C^{2}$, $f\in C(\mathcal{D})$ would have been more consistent with the story so far, but the $L^{\infty}$-assumption does not impede the upcoming analysis.

We briefly recall the arguments for existence and uniqueness of minimizers (the interested reader is referred to \cite[Chapter 8]{evans2010partial} for a lengthier discussion): The set $\mathcal{K}_{g,0}$ is closed in the strong (and hence also in the weak) $W^{1,2}$-topology.
Moreover, it is a convex subset of $W^{1,2}(\mathcal{D})$. $J$ is weakly lower semicontinuous, since the integrand is convex (and in fact also uniformly convex) in $\nabla u$, and $J$ is also coercive on $\mathcal{K}_{g,0}$.
Therefore there exists a unique minimizer in $\mathcal{K}_{g,0}$ to $J$.
Since the proof of uniqueness is so short we include it here: Let $u$ and $v$ be two minimizers and define $w\df\frac{u}{2}+\frac{v}{2}$ (which lies in $\mathcal{K}_{g,0}$ since this set is convex).
Then
\begin{align*}
J(w) & =\int_{\mathcal{D}}\left\lvert\nabla\left(\frac{u}{2}+\frac{v}{2}\right)\right\rvert^{2}+2\left(\frac{u}{2}+\frac{v}{2}\right)f\,dx\\
 & \le\frac{1}{2}\int_{\mathcal{D}}\lvert\nabla u\rvert^{2}+2uf\,dx+\frac{1}{2}\int_{\mathcal{D}}\lvert\nabla v\rvert^{2}+2vf\,dx\\
 & =\frac{J(u)}{2}+\frac{J(v)}{2}=\min_{\mathcal{K}_{g,0}}J.
\end{align*}
Since also $J(w)\ge\min_{\mathcal{K}_{g,0}}J$, equality holds in the inequality above, which by the strict convexity of $\xi\mapsto\lvert\xi\rvert^{2}$ happens if and only if $\nabla u=\nabla v$ a.e.~in $\mathcal{D}$, and since $u=v=g$ on the boundary we infer that $u=v$.

Even after reducing the obstacle to zero, we alter the functional once more so that the condition $u\ge0$ can be removed.
\begin{lem}
[{\cite[Lemma~1.4]{MR2962060}}]\label{lem:j-tilde} A function $u\in W^{1,2}(\mathcal{D})$ is a minimizer of $J$ in $\mathcal{K}_{g,0}$ if and only if $u$ is a minimizer to the functional \[ \tilde{J}(u)\df\int_{\mathcal{D}}\lvert\nabla u\rvert^{2}+2fu^{+}\,dx, \] in $\mathcal{K}_{g}\df\{u\in W^{1,2}(\mathcal{D}):u-g\in W_{0}^{1,2}(\mathcal{D})\}$.
\end{lem}

\begin{proof}
Let us start by noting that if $u\in\mathcal{K}_{g}$ then $u^{+}\in\mathcal{K}_{g,0}$ since the trace of $u^{-}$ is equal to the trace of $g^{-}$, which is zero by assumption.
Moreover, $\nabla u^{+}=\nabla u\cdot\chi_{\{u>0\}}$ by \cite[Lemma 7.6]{gilbarg:01}.
With this, it is straightforward to see that $\tilde{J}(u^{+})\le\tilde{J}(u)$ since the first term only decreases with this truncation, and the last term remains unchanged:
\begin{align*}
\tilde{J}(u) & =\int_{\mathcal{D}}\lvert\nabla u\rvert^{2}+2fu^{+}\,dx\ge\int_{\mathcal{D}}\lvert\nabla u\cdot\chi_{\{u>0\}}\rvert^{2}+2fu^{+}\,dx\\
 & =\int_{\mathcal{D}}\lvert\nabla u^{+}\rvert^{2}+2fu^{+}\,dx=\tilde{J}(u^{+}).
\end{align*}
Equality holds in the inequality above if and only if $\nabla u=0$ a.e.~on $\{u\le0\}$, which is equivalent to $\nabla u^{-}=0$ a.e.~in $\mathcal{D}$.
This in turn is equivalent to $u^{-}=0$ a.e.~in $\mathcal{D}$ since $u^{-}\in W_{0}^{1,2}(\mathcal{D})$ due to the fact that $g\ge0$ on $\partial\mathcal{D}$.

Therefore $\tilde{J}(u^{+})=\tilde{J}(u)$ if and only if $u\in\mathcal{K}_{g,0}$ (i.e.,~$u\ge0$ a.e.).
Moreover, $\tilde{J}$ and $J$ coincide on $\mathcal{K}_{g,0}$, hence also the minimizers coincide since the minimizer for $\tilde{J}$ does not lie in $\mathcal{K}_{g}\backslash \mathcal{K}_{g,0}$ by the previous discussion.
\end{proof}
We are now in a position to prove the local $W^{2,p}$-regularity of the minimizer, and the proof will require the standard tool of Calder\'{o}n--Zygmund estimates recalled below.
\begin{thm}
[Calder\'{o}n--Zygmund estimates, {\cite[Theorem 1.1]{MR2962060}}]\label{thm:cz-estimates} Let $f\in L^{p}(\mathcal{D})$, $p\in(1,\infty)$, and $u\in L^{1}(\mathcal{D})$ be a distributional solution to $\Delta u = f$ in $\mathcal{D}$.
Then $u\in W_{\text{loc}}^{2,p}(\mathcal{D})$ and \[ \lVert u\rVert_{W^{2,p}(K)}\le C\left(\lVert u\rVert_{L^{1}(\mathcal{D})}+\lVert f\rVert_{L^{p}(\mathcal{D})}\right) \] for all $K\Subset\mathcal{D}$, where $C=C(n,p,K,\mathcal{D})$.
\end{thm}

\begin{rem}
The Calder\'{o}n--Zygmund estimates above are \emph{not} true for $p=\infty$; otherwise, we would trivially obtain the desired $C_{\text{loc}}^{1,1}(\mathcal{D})$-regularity (or equivalently $W_{\text{loc}}^{2,\infty}(\mathcal{D})$-regularity).
\end{rem}

\begin{thm}
[$W^{2,p}$-regularity]\label{thm:w2p-regularity} Let $f\in L^\infty(\mathcal{D})$, and let $u$ be a solution to the obstacle problem in $\mathcal{D}$, i.e.,~a minimizer of $J$ in $\mathcal{K}_{g,0}$.
Then

\begin{enumerate}[label=(\roman{enumi})]
\item
$u\in W_{\text{loc}}^{2,p}(\mathcal{D})$, for any $p\in(1,\infty)$,
\item
$\Delta u=f\chi_{\{u>0\}}$ a.e.~in $\mathcal{D}$.
\end{enumerate}
\end{thm}

\begin{proof}
The proof of (i) is rather long so we divide it into steps.

\noindent\textbf{Step 1: Regularization.}
Consider the family of regularized problems,
\begin{equation}
\begin{cases}
\Delta u_{\epsilon}=f\chi_{\epsilon}(u_{\epsilon}) & \text{in }\mathcal{D},\\
\phantom{\Delta}u_{\epsilon}=g & \text{on }\partial \mathcal{D},
\end{cases}\tag{\ensuremath{*_{\epsilon}}}\label{eq:epsilon-reg}
\end{equation}
where $\chi_{\epsilon}$ is a smooth and nondecreasing approximation of the indicator function $\chi_{\{s>0\}}$ such that $\chi_{\{s>0\}}-\chi_{\epsilon}$ is supported on the interval $[-\epsilon,\epsilon]$.
By definition, a solution to \eqref{eq:epsilon-reg} is a critical point of the functional \[ J_{\epsilon}(v)\df\int_{\mathcal{D}}\lvert\nabla v\rvert^{2}+2f\Phi_{\epsilon}(v)\,dx \] in $\mathcal{K}_{g}$, where $\Phi_{\epsilon}(s)\df\int_{-\infty}^{s}\chi_{\epsilon}(t)\,dt$ is an approximation to the function $s_{+}$. $J_\epsilon$ is weakly lower semicontinuous and coercive since $\Phi_{\epsilon}(s) \le \lvert s \rvert + 1$, hence \[ J_\epsilon (v) \ge \lVert \nabla v \rVert_{L^2(\mathcal{D})}^2 - 2\lVert f\rVert_{L^{\infty}(\mathcal{D})}(\lVert v \rVert_{L^1(\mathcal{D})} + \lvert \mathcal{D}\rvert).
\] To finish the proof of coercivity we can use the Poincar\'{e} and Young inequalities, but the details are omitted since a similar argument will be used below to prove uniform $W^{1,2}$-estimates.
From weak lower semicontinuity and coercivity we infer the existence of a minimizer $u_{\epsilon}$ that satisfies
\begin{equation}
\int_{\mathcal{D}}\nabla u_{\epsilon}\cdot\nabla\eta+f\chi_{\epsilon}(u_{\epsilon})\eta\,dx=0\label{eq:variational}
\end{equation}
for all $\eta\in W_{0}^{1,2}(\mathcal{D})$.

\noindent\textbf{Step 2: Uniform $W^{1,2}$ bound.} We now claim that the
set of functions $\{u_{\epsilon}\}_{0<\epsilon<1}$ is uniformly bounded in $W^{1,2}(\mathcal{D})$, and in $W^{2,p}(K)$ for any $K\Subset\mathcal{D}$: Choose $\eta=u_{\epsilon}-g$ and plug into \eqref{eq:variational},
\begin{align}
0 & =\int_{\mathcal{D}}\nabla u_{\epsilon}\cdot\nabla(u_{\epsilon}-g)+f\chi_{\epsilon}(u_{\epsilon})(u_{\epsilon}-g)\,dx\nonumber \\
 & \ge\int_{\mathcal{D}}\lvert\nabla(u_{\epsilon}-g)\rvert^{2}-\underbrace{\lvert\nabla g\rvert\lvert\nabla(u_{\epsilon}-g)\rvert}_{\le\lvert\nabla g\rvert^{2}/2+\lvert\nabla(u_{\epsilon}-g)\rvert^{2}/2}-\lVert f\rVert_{L^{\infty}(\mathcal{D})}\lvert u_{\epsilon}-g\rvert\,dx\nonumber \\
 & \ge\int_{\mathcal{D}}\frac{\lvert\nabla(u_{\epsilon}-g)\rvert^{2}}{2}-\frac{\lvert\nabla g\rvert^{2}}{2}-\lVert f\rVert_{L^{\infty}(\mathcal{D})}\lvert u_{\epsilon}-g\rvert\,dx.\label{eq:use-young}
\end{align}
In the first inequality, we used that $\lvert\chi_{\epsilon}(u_{\epsilon})\rvert\le1$.
For the last term, we employ the Young inequality $ab\le \tau a^{2} + b^{2}/(4\tau)$ with $\tau\df C_P^2$ where $C_{P}$ is the constant in the Poincar\'{e} inequality $\int_{\mathcal{D}}\lvert v\rvert^{2}\,dx\le C_{P}^{2}\int_{\mathcal{D}}\lvert\nabla v\rvert^{2}\,dx$ for $v\in W_{0}^{1,2}(\mathcal{D})$.
Consequently the expression in \eqref{eq:use-young} can be estimated from below by
\begin{align*}
 & \int_{\mathcal{D}}\frac{\lvert\nabla(u_{\epsilon}-g)\rvert^{2}}{2}-\frac{\lvert\nabla g\rvert^{2}}{2}-C_{P}^{2}\lVert f\rVert_{L^{\infty}(\mathcal{D})}^{2}-\frac{\lvert u_{\epsilon}-g\rvert^{2}}{4C_{P}^{2}}\,dx\\
\ge & \{\text{Poincar\'{e}'s inequality}\}\\
\ge & \int_{\mathcal{D}}\frac{\lvert\nabla(u_{\epsilon}-g)\rvert^{2}}{4}-\frac{\lvert\nabla g\rvert^{2}}{2}-C_{P}^{2}\lVert f\rVert_{L^{\infty}(\mathcal{D})}^{2}\,dx.
\end{align*}
Since we know that the last expression is bounded from above by zero, we conclude that \[ \lVert\nabla(u_{\epsilon}-g)\rVert_{L^{2}(\mathcal{D})}^2\le2\lVert\nabla g\rVert_{L^{2}(\mathcal{D})}^{2}+4C_{P}^{2}\lvert\mathcal{D}\rvert\lVert f\rVert_{L^{\infty}(\mathcal{D})}^{2} \] from which it is easily inferred that $\lVert u_{\epsilon}\rVert_{W^{1,2}(\mathcal{D})}$ is uniformly bounded by a quantity that only depends on $f$, $g$, and $\mathcal{D}$.
Therefore there exists a sequence $\epsilon_{k}\to0$ as $k\to\infty$, and $w\in W^{1,2}(\mathcal{D})$ such that
\begin{align*}
u_{\epsilon_{k}} & \rightharpoonup w,\phantom{(x)}\qquad\text{in }W^{1,2}(\mathcal{D}),\\
u_{\epsilon_{k}} & \to w,\phantom{(x)}\qquad\text{in }L^{2}(\mathcal{D}),\\
u_{\epsilon_{k}}(x) & \to w(x),\qquad\text{for a.e.~\ensuremath{x} in }\mathcal{D}.
\end{align*}
The function $w$ lies in $\mathcal{K}_{g}$ since $\mathcal{K}_{g}$ is weakly closed.

\noindent\textbf{Step 3: Compactness and $W_\text{loc}^{2,p}$ convergence.}
Now, given $K\Subset\mathcal{D}$, by an application of the Calder\'{o}n--Zygmund estimates in Theorem~\ref{thm:cz-estimates} \[ \lVert u_{\epsilon}\rVert_{W^{2,p}(K)}\le C(\lVert u_{\epsilon}\rVert_{L^{1}(\mathcal{D})}+\lVert f\rVert_{L^{\infty}(\mathcal{D})})\le C(K,\mathcal{D},f,g).
\] Therefore there exists a further subsequence of $\epsilon_{k}$ (which we continue to denote by $\epsilon_{k}$ to avoid cumbersome notation) such that, for all $p\in(1,\infty)$, $u_{\epsilon_{k}}\rightharpoonup w$ in $W_{\text{loc}}^{2,p}(\mathcal{D})$.
Hence $w\in W_{\text{loc}}^{2,p}(\mathcal{D})$ for any $p\in(1,\infty)$, and due to lower semicontinuity with respect to weak convergence, and also the dominated convergence theorem,
\begin{align*}
\int_{\mathcal{D}}\lvert\nabla w\rvert^{2}\,dx & \le\liminf_{k\to\infty}\int_{\mathcal{D}}\lvert\nabla u_{\epsilon_{k}}\rvert^{2}\,dx,\\
\int_{\mathcal{D}}f(x)w^{+}(x)\,dx & =\lim_{k\to\infty}\int_{\mathcal{D}}f\Phi_{\epsilon_{k}}(u_{\epsilon_{k}})\,dx.
\end{align*}
These facts combined show that $w$ minimizes the functional $\tilde{J}$, \[ \tilde{J}(w)\le\liminf_{\epsilon_{k}\to0}J_{\epsilon_{k}}(u_{\epsilon_{k}})\le\liminf_{\epsilon_{k}\to0}J_{\epsilon_{k}}(v)\le\tilde{J}(v),\qquad\forall v\in\mathcal{K}_{g}.
\] By Lemma~\ref{lem:j-tilde}, $w$ is a solution to the obstacle problem, and has to coincide with $u$ due to uniqueness.

To show (ii), we use that $u\in W_{\text{loc}}^{2,p}(\mathcal{D})$: On the coincidence set we have $\nabla u=0$ and $D^{2}u=0$ a.e.
Consequently, $\Delta u=\trace(D^{2}u) = 0$ a.e.~on $\{x\in\mathcal{D}:u=0\}$.
Now let us consider the set $\{x\in\mathcal{D}:u>0\}$.
Because of the embedding of Sobolev spaces into H\"{o}lder spaces for $p>n$, $u\in C_{\text{loc}}^{1,\alpha}(\mathcal{D})$ for any $\alpha\in(0,1)$.
Therefore the set $\{x\in\mathcal{D}:u(x)>0\}$ is open; take a ball $B\Subset\{x\in\mathcal{D}:u(x)>0\}$, $\eta\in C_{c}^{\infty}(B)$ and $\epsilon>0$ small enough that $v_{\epsilon}\df u\pm\epsilon\eta>0$ in $B$.
Since $\tilde{J}(u)\le\tilde{J}(v_{\epsilon})$ and
\begin{align*}
\tilde{J}(u) & =\int_{\mathcal{D}}\lvert\nabla u\rvert^{2}+2fu^{+}\,dx,\\
\tilde{J}(v_{\epsilon}) & =\int_{\mathcal{D}}\lvert\nabla u\rvert^{2}\pm2\epsilon\nabla u\cdot\nabla\eta+\epsilon^{2}\lvert\nabla\eta\rvert^{2}+2f(u\pm\epsilon\eta)^{+}\,dx,
\end{align*}
we deduce that
\begin{equation}
\int_{\mathcal{D}}fu^{+}\,dx\le\int_{\mathcal{D}}f(u\pm\epsilon\eta)^{+}\pm\epsilon\nabla u\cdot\nabla\eta+\frac{\epsilon^{2}\lvert\nabla\eta\rvert^{2}}{2}\,dx.\label{eq:eta-epsilon-variation}
\end{equation}
Since $\eta\in C_{c}^{\infty}(B)$, $u>0$ in $B$ and $u\pm\epsilon\eta>0$ in $B$,
\begin{align*}
\int_{\mathcal{D}}f(u\pm\epsilon\eta)^{+} & =\int_{\mathcal{D}\backslash B}fu^{+}+\int_{B}f(u\pm\epsilon\eta)\,dx,\\
\int_{\mathcal{D}}fu^{+}\,dx & =\int_{\mathcal{D}\backslash B}fu^{+}\,dx+\int_{B}fu\,dx.
\end{align*}
We plug these expressions into \eqref{eq:eta-epsilon-variation} to get \[ 0\le\int_{B}\pm\epsilon(f\eta+\nabla u\cdot\nabla\eta)+\frac{\epsilon^{2}\lvert\nabla\eta\rvert^{2}}{2}\,dx.
\] Divide by $\epsilon$, let $\epsilon\to0^{+}$ and integrate by parts, \[ 0\le\pm\int_{B}\eta(f-\Delta u)\,dx, \] which implies that $\Delta u=f$ a.e.~on $\{x\in\mathcal{D}:u(x)>0\}$.
\end{proof}

In view of the preceding theorem, every minimizer satisfies \[
\begin{cases}
	\Delta u=f\chi_{\{u>0\}},\\
	\phantom{\Delta}u\ge0,
\end{cases}
\quad \text{a.e.~in } \mathcal{D}.
\] Since the results that follow are local, we drop the boundary condition and henceforth work with this PDE formulation: by a \emph{solution to the obstacle problem} we mean a function $u\in W_{\text{loc}}^{2,p}(\mathcal{D})$, $p\in(1,\infty)$, satisfying the system above.
Nothing is lost in this passage when $f\ge0$, since a solution with boundary values $g$ is then also the minimizer.
For sign-changing $f$, however, the PDE formulation is strictly weaker, as it only sees the values of $f$ in $\{u>0\}$ (see Remark~\ref{rem:degenerate}).

\subsection{Optimal \texorpdfstring{$C^{1,1}$}{C1,1}-Regularity}

In the previous section, we established that the solution $u$ to the obstacle problem belongs not only to $W^{1,2}(\mathcal{D})$, but also to $W_{\text{loc}}^{2,p}(\mathcal{D})$ for any $p\in(1,\infty)$.
Hence it is also $C_{\text{loc}}^{1,\alpha}(\mathcal{D})$ for any $\alpha\in(0,1)$.
However, as noted in Section~\ref{sec:main-section}, the derivatives of $u$ are actually Lipschitz, suggesting that under suitable assumptions we can reach the optimal regularity \[ C_{\text{loc}}^{1,1}(\mathcal{D})=W_{\text{loc}}^{2,\infty}(\mathcal{D}) \] (see \cite[Chapter 5.8]{evans2010partial} for a proof of this equivalence).

To begin, we prove a sub-quadratic growth lemma.
It states that the function $u$ cannot grow faster than quadratically away from the free boundary $\Gamma(u)=\partial\Omega\cap\mathcal{D}$, where we recall the noncoincidence set $\Omega(u)\df\{x\in\mathcal{D}:u(x)>0\}$ and the coincidence set $\Lambda(u)\df\{x\in\mathcal{D}:u(x)=0\}$.
\begin{lem}
[Sub-quadratic growth {\cite[Theorem 2.1]{MR2962060}}]\label{lem:sub-quadratic-growth} Let $f\in L^{\infty}(\mathcal{D})$ and $u$ be a solution to the obstacle problem \[
\begin{cases}
	\Delta u=f\chi_{\{u>0\}},\\
	\phantom{\Delta}u\ge0,
\end{cases}
\quad \text{in } \mathcal{D}.
\] Let $x^{0}\in\Gamma(u)$ and $B_{2R}(x^{0})\subset\mathcal{D}$.
Then \[ \sup_{B_{R}(x^{0})}u\le C\lVert f\rVert_{L^{\infty}(\mathcal{D})}R^{2}.
\]
\end{lem}

\begin{proof}
Without loss of generality, assume $x^{0}=0$.
Decompose $u=u_{1}+u_{2}$ in $B_{2R}$, where $u_{1}$ and $u_{2}$ are defined as follows:
\begin{align*}
\begin{cases}
\Delta u_{1}=\Delta u,\quad\Delta u_{2}=0 & \text{in }B_{2R},\\
\phantom{\Delta}u_{1}=0,\qquad\phantom{\Delta}u_{2}=u & \text{on }\partial B_{2R}.
\end{cases}
\end{align*}
We also define $\phi\df\frac{1}{2n}(4R^{2}-\lvert x\rvert^{2})$, which then solves \[
\begin{cases}
\Delta\phi=-1 & \text{in }B_{2R},\\
\phantom{\Delta}\phi=0 & \text{on }\partial B_{2R}.
\end{cases}
\] By the comparison principle \[ -\lVert f\rVert_{L^{\infty}(\mathcal{D})}\phi(x)\le u_{1}(x)\le\lVert f\rVert_{L^{\infty}(\mathcal{D})}\phi(x) \] in $B_{2R}$, which implies that $\lvert u_{1}\rvert\le C\lVert f\rVert_{L^{\infty}(\mathcal{D})}R^{2}$.
Since $u_2$ is harmonic in $B_{2R}$ and nonnegative on the boundary, the maximum principle gives $u_2 \ge 0$ in $B_{2R}$.
Moreover, $u(0) = 0$ implies $u_2(0)=-u_1(0)\le C\lVert f\rVert_{L^{\infty}(\mathcal{D})}R^{2}$.
By the Harnack inequality, \[ \sup_{B_{R}}u_{2}\le Cu_{2}(0)\le C\lVert f\rVert_{L^{\infty}(\mathcal{D})}R^{2}.
\] Combining these, \[ \sup_{B_{R}}u\le\sup_{B_{R}}u_{1}+\sup_{B_{R}}u_{2}\le C\lVert f\rVert_{L^{\infty}(\mathcal{D})}R^{2}.
\]
\end{proof}
At this stage, achieving $C^{1,1}$-regularity requires more than $f\in L^{\infty}$.
In fact, continuity of $f$ alone is insufficient.
The sharp condition is that $f$ admits a so-called $C^{1,1}$-potential as defined below.
\begin{defn}
[$C^{1,1}$-potential] A function $f\in L^{\infty}(\mathcal{D})$ has a $C^{1,1}$-potential in $\mathcal{D}$ if there exists a function $\psi\in C^{1,1}(\mathcal{D})$ such that $\Delta\psi=f$ in $\mathcal{D}$.
\end{defn}

\begin{rem}
We reuse the symbol $\psi$: it originally denoted the obstacle, for which
$f=-\Delta\psi$, whereas here it denotes a potential satisfying $\Delta\psi=f$. 
The two differ by a sign, and since only the $C^{1,1}$-regularity will enter the analysis, requiring $f$ to have a $C^{1,1}$-potential is thus equivalent to requiring the obstacle to be $C^{1,1}$.
\end{rem}

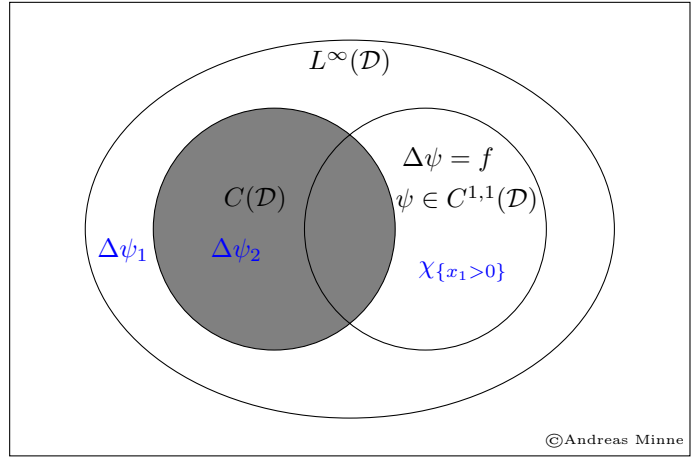
\begin{figure}
\begin{centering}
\begin{tikzpicture}[fill=gray]
\fill (0.5,0) circle (1.6);
\draw (1.5,0) ellipse (3.5cm and 2.5cm) (1.5,2.5) node [text=black,below] {$L^\infty(\mathcal{D})$};
\draw (0.5,0) circle (1.6) (0.25,0.7)  node [text=black,below] {$C(\mathcal{D})$}       
(2.5,0) circle (1.6) (2.8,1.2)  node [text=black,below] {$\Delta \psi = f$}
(3.05,0.7) node [text=black,below] {$\psi \in C^{1,1}(\mathcal{D})$}      
(-3,-3) rectangle (6,3)
(-1.5,0) node [text=blue,below] {$\Delta \psi_1$}
(0,0) node [text=blue,below] {$\Delta\psi_2$}
(3,-0.3) node [text=blue,below] {$\chi_{\{x_1>0\}}$};

\node at (5,-2.8) (copyright-text) {\tiny{\copyright Andreas Minne}};

\end{tikzpicture}\caption{Venn diagram illustrating the relation between $L^{\infty}$-functions,
continuous functions, and functions with a $C^{1,1}$-potential.
The example functions $\psi_{1}$, $\psi_{2}$, and the characteristic function (all written in blue) are placed in their corresponding function space.
Note that any other $C^{1,1}$-potential of $\Delta \psi_{1}$, and $\Delta \psi_{2}$ differ by most a harmonic function, which is smooth, \label{fig:venn-diagram} } \par\end{centering}
\end{figure}

To emphasize the distinction between functions which have a $C^{1,1}$ potential, and functions which do not, we give three examples: one where $f$ has a $C^{1,1}$-potential but is not continuous, one where $\Delta\psi=f\in L^{\infty}(\mathcal{D})$ but $\psi\not\in C^{1,1}(\mathcal{D})$, and one where $\Delta\psi=f\in C(\mathcal{D})$ but $\psi\not\in C^{1,1}(\mathcal{D})$.
The relation between $L^{\infty}$-functions, continuous functions, and functions with a $C^{1,1}$-potential is illustrated in the Venn diagram in Figure~\ref{fig:venn-diagram}, where the $L^{\infty}$-functions are represented by the ellipse, the continuous functions by the gray disk on the left, and the functions with a $C^{1,1}$-potential by the right disk.
The example functions are included (in blue) as well in the Venn diagram. 

\begin{example}
[{$f$ has a $C^{1,1}$-potential but $f\not\in C(\mathcal{D})$}] Let $f=\chi_{\{x_{1}>0\}}$ where $x_{1}$ is the first coordinate.
Then $f$ has a jump discontinuity across the hyperplane $\{x\in\mathbb{R}^{n}:x_{1}=0\}$.
However, $\psi=(x_{1}^{+})^{2}/2$ is a $C^{1,1}$-potential for $f$.
\end{example}

\begin{example}
[{$f\in L^\infty(\mathcal{D})$ but $\psi\not\in C^{1,1}(\mathcal{D})$}] Consider $\psi_{1}(x,y)\df xy\ln(x^{2}+y^{2})$ for $x$, $y\in\mathbb{R}$.
Then \[ f\df\Delta\psi_{1}=\frac{8xy}{x^{2}+y^{2}}\in L^{\infty}(B_{1/2})\backslash C(B_{1/2}) \] while the mixed partial derivatives \[ \partial_{xy}\psi_{1}=\underbrace{2\frac{(x^{4}+y^{4})}{(x^{2}+y^{2})^{2}}}_{\text{bounded}}+\ln(x^{2}+y^{2})\to-\infty \] as we approach the origin.
This shows that $\psi_{1}\not\in W_{\text{loc}}^{2,\infty}(B_{1/2})=C_{\text{loc}}^{1,1}(B_{1/2})$.
\end{example}

\begin{example}
[{$f\in C(\mathcal{D})$ but $\psi\not\in C^{1,1}(\mathcal{D})$}] Consider $\psi_{2}(x,y)\df xy(-\ln(x^{2}+y^{2}))^{1/2}$ for $x$, $y\in\mathbb{R}$.
Then 
\[ f\df\Delta\psi_{2}=\frac{xy(4\ln(x^{2}+y^{2})-1)}{(x^{2}+y^{2})(-\ln(x^{2}+y^{2}))^{3/2}}\in C(B_{1/2}) 
\] 
(we have naturally defined $f$ to be zero at the origin) while the mixed partial derivatives 
\[ \partial_{xy}\psi_{2}=\underbrace{\frac{(x^{4}+y^{4})\ln(x^{2}+y^{2})-x^{2}y^{2}}{(x^{2}+y^{2})^{2}(-\ln(x^{2}+y^{2}))^{3/2}}}_{\text{bounded}}+(-\ln(x^{2}+y^{2}))^{1/2}\to \infty 
\] 
as we approach the origin.
This shows that $\psi_{2}\not\in W_{\text{loc}}^{2,\infty}(B_{1/2})=C_{\text{loc}}^{1,1}(B_{1/2})$.
\end{example}

Before we turn to the proof of $C^{1,1}$-regularity we need the following lemma which relates two functions for which the Laplacians coincide on a ball.
\begin{lem}
\label{lem:harmonic-relation} If $\Delta v=\Delta\psi$ in $B_{2R}(x^{0})$
and $\psi\in C^{1,1}(B_{2R}(x^{0}))$, then \[ \lVert D^{2}v\rVert_{L^{\infty}(B_{R}(x^{0}))}\le C\bigg(\frac{\lVert v\rVert_{L^{\infty}(B_{2R}(x^{0}))}}{R^{2}}+\lVert D^{2}\psi\rVert_{L^{\infty}(B_{2R}(x^{0}))}\bigg).
\]
\end{lem}

\begin{proof}
To ease the notation, let $B_{R}=B_{R}(x^{0})$ and $B_{2R}=B_{2R}(x^{0})$, and consider the difference $w:=v-\psi$ which is harmonic.
Then for $w$ (see Chapter 2 in \cite{gilbarg:01}) we have
\begin{equation}
\lVert D^{2}w\rVert_{L^{\infty}(B_{R})}\le C\frac{\lVert w\rVert_{L^{\infty}(B_{2R})}}{R^{2}}\le C\bigg(\frac{\lVert v\rVert_{L^{\infty}(B_{2R})}}{R^{2}}+\frac{\lVert\psi\rVert_{L^{\infty}(B_{2R})}}{R^{2}}\bigg).\label{eq:harmonicestimate}
\end{equation}
We can without loss of generality assume that $\psi(x^{0})=\lvert\nabla\psi(x^{0})\rvert=0$ (otherwise redefine $\psi$ as $\psi(x)-\psi(x^{0})-\nabla\psi(x^{0})\cdot(x-x^{0})$).
Then the mean value inequality gives, for $x\in B_{2R}$, \[ \lvert\psi(x)\rvert=\lvert\psi(x)-\psi(x^{0})\rvert\le\lVert \nabla\psi\rVert_{L^{\infty}(B_{2R})}\lvert x-x^{0}\rvert\le\lVert \nabla\psi\rVert_{L^{\infty}(B_{2R})}2R, \] and similarly for the gradient, \[ \lvert \nabla\psi(x)\rvert=\lvert \nabla\psi(x)-\nabla\psi(x^{0})\rvert\le\lVert D^{2}\psi\rVert_{L^{\infty}(B_{2R})}\lvert x-x^{0}\rvert\le\lVert D^{2}\psi\rVert_{L^{\infty}(B_{2R})}2R. \] Therefore $\lVert\psi\rVert_{L^{\infty}(B_{2R})}\le4R^{2}\lVert D^{2}\psi\rVert_{L^{\infty}(B_{2R})}$ which we put in \eqref{eq:harmonicestimate},
\begin{align*}
\lVert D^{2}w\rVert_{L^{\infty}(B_{R})} & \le C\bigg(\frac{\lVert v\rVert_{L^{\infty}(B_{2R})}}{R^{2}}+\frac{\lVert\psi\rVert_{L^{\infty}(B_{2R})}}{R^{2}}\bigg)\\
 & \le C\bigg(\frac{\lVert v\rVert_{L^{\infty}(B_{2R})}}{R^{2}}+\lVert D^{2}\psi\rVert_{L^{\infty}(B_{2R})}\bigg).
\end{align*}
The same inequality holds for $D^2 v$ in place of $D^2 w$ if we add an $L^\infty$-term of $D^2 \psi$ on the right hand side.
\end{proof}
We now arrive at the final step of our regularity analysis, namely the optimal local $C^{1,1}$-regularity (or, equivalently, $W^{2,\infty}$-regularity) of solutions to the obstacle problem.
\begin{thm}
[Optimal $C^{1,1}$-regularity {\cite[Theorem 2.3]{MR2962060}}]\label{thm:c11-regularity} Let $f$ have a $C^{1,1}$-potential $\psi$, and assume $u$ to be a solution to the obstacle problem \[
\begin{cases}
	\Delta u=f\chi_{\{u>0\}},\\
	\phantom{\Delta}u\ge0,
\end{cases}
\quad \text{in } \mathcal{D}.
\] Then $u\in C_{\text{loc}}^{1,1}(\mathcal{D})$ and, for $K\Subset\mathcal{D}$, \[ \lVert u\rVert_{C^{1,1}(K)}\le C(\lVert u\rVert_{L^{\infty}(\mathcal{D})}+\lVert D^{2}\psi\rVert_{L^{\infty}(\mathcal{D})}).
\]
\end{thm}

\begin{proof}
If $\Lambda(u)$ is empty, $\Delta u = f$ throughout $\mathcal{D}$, so $u$ is $C_\text{loc}^{1,1}(\mathcal{D})$ by Lemma~\ref{lem:harmonic-relation} since $f$ has a $C^{1,1}$-potential. Hence we will assume that $\Lambda(u)\ne \varnothing$. 
Since $u\in W_{\text{loc}}^{2,p}(\mathcal{D})$ and $D^{2}u=0$ a.e.~on the coincidence set $\Lambda(u)$, we consider $x^{0}\in\Omega(u)\cap K$ for $K\Subset\mathcal{D}$.
Let $d\df\dist(x^{0},\Lambda(u))$, $\delta\df\dist(K,\partial\mathcal{D})$ (note that $\delta,$ $d\ne \infty$ since $\mathcal{D}$ is bounded and  $\Lambda(u)$ is non-empty).
Then we get the following cases:

\noindent\textbf{Case (i)} ($d<\delta/5$). Let $y^{0}\in\partial B_{d}(x^{0})\cap\Gamma(u)$.
Then $B_{4d}(y^{0})\subset B_{5d}(x^{0})\Subset \mathcal{D}$.
The sub-quadratic growth, Lemma~\ref{lem:sub-quadratic-growth}, gives
\begin{equation}
\lVert u\rVert_{L^{\infty}(B_{2d}(y^{0}))}\le C\lVert f\rVert_{L^{\infty}(\mathcal{D})}(2d)^{2}.\label{eq:sub-quadratic-c11}
\end{equation}
Since $B_{d}(x^{0})$ lies completely in the noncoincidence set by definition, $\Delta u=f=\Delta\psi$ in $B_{d}(x^{0})$, so we can apply Lemma~\ref{lem:harmonic-relation} and \eqref{eq:sub-quadratic-c11},
\begin{align*}
\lVert D^{2}u\rVert_{L^{\infty}(B_{d/2}(x^{0}))} & \le C\bigg(\frac{\lVert u\rVert_{L^{\infty}(B_{d}(x^{0}))}}{d^{2}}+\lVert D^{2}\psi\rVert_{L^{\infty}(B_{d}(x^{0}))}\bigg)\\
 & \le C\bigg(\lVert f\rVert_{L^{\infty}(\mathcal{D})}+\lVert D^{2}\psi\rVert_{L^{\infty}(B_{d}(x^{0}))}\bigg)\\
 & \le C\lVert D^{2}\psi\rVert_{L^{\infty}(\mathcal{D})}
\end{align*}
where we have $\lVert f\rVert_{L^{\infty}(\mathcal{D})}\le C\lVert D^{2}\psi\rVert_{L^{\infty}(\mathcal{D})}$ since $f$ is just $\trace(D^{2}\psi)$.

\noindent\textbf{Case (ii)} ($d\ge\delta/5$). Set $\rho\df\min(d,\delta/2)$. Then $B_{\rho}(x^{0})\subseteq\Omega(u)$ since $\rho\le d$, and $\overline{B_{\rho}(x^{0})}\Subset\mathcal{D}$ since $\rho\le\delta/2$. Hence $\Delta u=f$ in $B_{\rho}(x^{0})$, and Lemma~\ref{lem:harmonic-relation} gives
\begin{align*}
	\lVert D^{2}u\rVert_{L^{\infty}(B_{\rho/2}(x^{0}))} & \le C\bigg(\frac{\lVert u\rVert_{L^{\infty}(B_{\rho}(x^{0}))}}{\rho^{2}}+\lVert D^{2}\psi\rVert_{L^{\infty}(B_{\rho}(x^{0}))}\bigg)\\
	& \le C\bigg(\frac{\lVert u\rVert_{L^{\infty}(\mathcal{D})}}{\delta^{2}}+\lVert D^{2}\psi\rVert_{L^{\infty}(\mathcal{D})}\bigg),
\end{align*}
where the last inequality uses $\rho\ge\delta/5$.

From the combination of these two cases, classical $C^{1,\alpha}$ estimates of $u$, and the fact that $\delta$ only depends on $K$, 
\[ 
\lVert u\rVert_{C^{1,1}(K)}\le C\Big(\lVert u\rVert_{L^{\infty}(\mathcal{D})}+\lVert D^{2}\psi\rVert_{L^{\infty}(\mathcal{D})}\Big).
\]
\end{proof}

\section{Analysis of the Free Boundary}\label{sec:free-boundary}

Now that we have established the optimal regularity of the solution, we turn our attention to the free boundary.
Characterizing the free boundary is, however, a highly nontrivial task, and our analysis will necessarily be limited in scope.
The interested reader is again referred to the book by Petrosyan et al.~\cite{MR2962060} for a more complete analysis, and there has also been some recent progress by Figalli et al.~\cite{MR3904453,FR-OS20}.

\begin{figure}
\begin{minipage}[t]{0.45\textwidth}%
\centering
\begin{tikzpicture}
\node[anchor=south east,inner sep=0] (image) at (0,0) {\includegraphics[width=\textwidth]{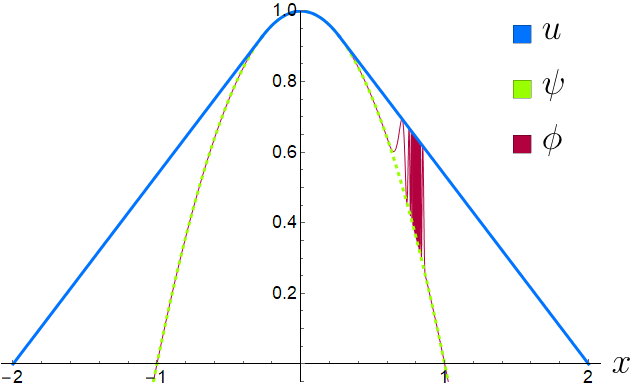}};
\node[below left] (copyright-text) {\tiny{\copyright Andreas Minne}};
\end{tikzpicture}\subcaption{}\label{fig:mean-obstacle}%
\end{minipage}\hfill{}%
\begin{minipage}[t]{0.45\textwidth}%
\centering
\begin{tikzpicture}
\node[anchor=south east,inner sep=0] (image) at (0,0) {\includegraphics[width=\textwidth]{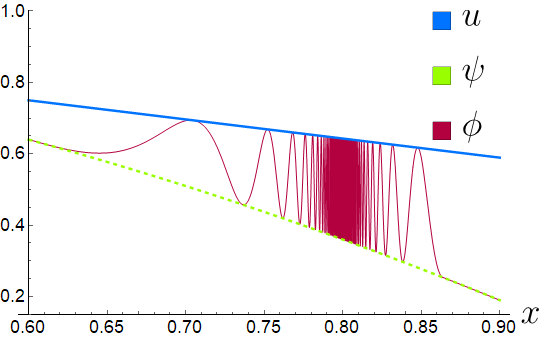}};
\node[below left] (copyright-text) {\tiny{\copyright Andreas Minne}};
\end{tikzpicture}\subcaption{}\label{fig:mean-obstacle-zoom}%
\end{minipage}

\caption{A troublesome obstacle $\phi$ that renders $\Gamma(u)$ rather messy.}
\end{figure}

Just to illustrate how things can go sour, let us return to the very first example with the rubber band, Figure~\ref{fig:1Dslnwithobstacle}.
In that example there are exactly two free boundary points, at $x=\pm(2-\sqrt{3})$.
It is, however, possible to (nota bene, smoothly!) alter the obstacle $\psi$ to another obstacle $\phi$ in such a way that the solution remains unchanged, while the free boundary becomes pathological (for instance, \emph{any} prescribed closed set).
A suggestive sketch of this idea is found in Figure~\ref{fig:mean-obstacle} (and \ref{fig:mean-obstacle-zoom} for a zoomed-in version), and we can only imagine the horrors that may occur in higher dimensions.

To avoid this issue we continue our analysis under the assumption that $f$ is uniformly bounded away from zero, \[ f\ge c>0.
\] Recall that this assumption is equivalent to assuming that the obstacle $\psi$ is strictly superharmonic, $\Delta\psi=-f\le-c<0$, in the original formulation of the problem.
With this assumption, we can improve the sub-quadratic growth in Lemma~\ref{lem:sub-quadratic-growth} to quadratic growth in the noncoincidence set and on the free boundary.
The super-quadratic property below is often referred to as nondegeneracy.
\begin{lem}
[Nondegeneracy {\cite[Lemma 3.1]{MR2962060}}]\label{lem:nondegeneracy} Let $f\in L^{\infty}(\mathcal{D})$, $f\ge c>0,$ and assume $u$ to be a solution to the obstacle problem \[
\begin{cases}
	\Delta u=f\chi_{\{u>0\}},\\
	\phantom{\Delta}u\ge0,
\end{cases}
\quad \text{in } \mathcal{D}.
\] If $x^{0}\in\overline{\Omega(u)}$, then \[ \sup_{B_{r}(x^{0})}u\ge u(x^{0})+\frac{c}{2n}r^{2},\qquad B_{r}(x^{0})\Subset\mathcal{D}.
\]
\end{lem}

\begin{proof}
Let $x^{0}\in\Omega(u)$ and consider \[ w(x)\df u(x)-u(x^{0})-\frac{c\lvert x-x^{0}\rvert^{2}}{2n}.
\] In $B_{r}(x^{0})\cap\Omega(u)$, we have $\Delta w\ge0$, and also $w(x^{0})=0$.
The maximum principle implies that \[ 0\le\sup_{B_{r}(x^{0})\cap\Omega(u)}w=\max_{\partial(B_{r}(x^{0})\cap\Omega(u))}w.
\] We have the set inclusion \[ \partial(B_{r}(x^{0})\cap\Omega(u))\subseteq(\partial B_{r}(x^{0})\cap\overline{\Omega(u)})\cup(\overline{B_{r}(x^{0})}\cap\partial\Omega(u)), \] and from the fact that $w<0$ on $\overline{B_{r}(x^{0})}\cap\partial\Omega(u)$ we infer that \[ 0\le\max_{\partial B_{r}(x^{0})\cap\Omega(u)}w=\max_{\partial B_{r}(x^{0})\cap\Omega(u)}u-u(x^{0})-\frac{c}{2n}r^{2}, \] hence \[ \sup_{B_{r}(x^{0})}u\ge u(x^{0})+\frac{c}{2n}r^{2}.
\] For $x^{0}\in\Gamma(u)$, let $\Omega(u)\ni x^{k}\to x^{0}$ as $k\to\infty$, and use continuity of $u$.
\end{proof}
\begin{rem}\label{rem:degenerate}
An easy example of a degenerate solution, where the condition $f\ge c>0$ is violated, is $(x^{+})^{3}/6$, which solves \[
\begin{cases}
\Delta u=\lvert x\rvert\chi_{\{u>0\}} & \text{in }(-1,1),\\
\phantom{\Delta}u\ge0 & \text{in }(-1,1),\\
u(-1)=0, & u(1)=1/6.
\end{cases}
\]
Here $f=\lvert x\rvert$ vanishes at the free boundary point $x=0$, and the growth is cubic rather than quadratic.
Note that $u$ also solves the same system with $f=x$ in place of $\lvert x\rvert$, since $f$ is only seen in $\{u>0\}$; but then $u$ is \emph{not} a minimizer, since $f<0$ on the coincidence set $(-1,0]$ and perturbing $u$ upwards there decreases the energy.
\end{rem}

As a consequence of Lemma~\ref{lem:nondegeneracy}, we also get a linear nondegeneracy of the gradient.
\begin{cor}
[Nondegeneracy of the gradient {\cite[Corollary 3.7]{MR2962060}}]\label{cor:nondegeneracy-gradient} Let $f\in L^{\infty}(\mathcal{D})$, $f\ge c>0$, and assume $u$ to be a solution to the obstacle problem \[
\begin{cases}
\Delta u=f\chi_{\{u>0\}},\\
\phantom{\Delta}u\ge0,
\end{cases}
\] in $\mathcal{D}$.
If $x^{0}\in\overline{\Omega(u)}$, then \[ \sup_{B_{r}(x^{0})}\lvert\nabla u\rvert\ge\frac{c}{2n}r,\qquad B_{r}(x^{0})\Subset\mathcal{D}.
\]
\end{cor}

\begin{proof}
The gradient $\nabla u$ is continuous in $\mathcal{D}.$ By nondegeneracy of $u$, Lemma~\ref{lem:nondegeneracy}, choose a point $x\in\partial B_{r}(x^{0})$ such that $\frac{c}{2n}r^{2}\le u(x)-u(x^{0})$.
Then \[ \frac{c}{2n}r^{2}\le u(x)-u(x^{0})\le\sup_{B_{r}(x^{0})}\lvert\nabla u\rvert\lvert x-x^{0}\rvert=\sup_{B_{r}(x^{0})}\lvert\nabla u\rvert r \] from which we get the desired inequality by dividing by $r$ on both sides.
\end{proof}
The first result concerning the free boundary is how big it is; we would like to prove that it indeed is an interface between the coincidence and noncoincidence set in the sense that it does not have full measure.
To do this we prove a stronger property --- porosity --- which we define below.
\begin{defn}
A measurable set $E\subseteq\mathbb{R}^{n}$ is \emph{$\delta$-porous}, $0<\delta<1$ if every ball $B_{r}(x^{0})\subset\mathbb{R}^{n}$ contains a smaller ball $B_{\delta r}(y)\subset B_{r}(x^{0})\backslash E$.
Moreover, $E$ is \emph{locally porous} in an open set $\mathcal{D}\subseteq\mathbb{R}^{n}$ if $E\cap K$ is $\delta$-porous for any $K\Subset\mathcal{D}$, with $\delta=\delta(K)$.
\end{defn}

\begin{lem}
[Porosity implies zero measure {\cite[Proposition 3.9]{MR2962060}}]If $E\subset\mathbb{R}^{n}$ is locally porous in $\mathcal{D}\subseteq\mathbb{R}^{n}$, then $\lvert E\cap\mathcal{D}\rvert=0$.
\end{lem}

\begin{proof}
Let $x^{0}\in E$ and $B_{r}(x^{0})\subset K\Subset\mathcal{D}$.
Since $B_{r}(x^{0})$ is the disjoint union of $B_{r}(x^{0})\cap E$ and $B_{r}(x^{0})\backslash E$, and the latter at least contains a ball $B_{\delta r}(y)$ of size $(\delta r)^{n}$, we have
\begin{align*}
\lvert E\cap B_{r}(x^{0})\rvert & =\lvert B_{r}(x^{0})\rvert-\lvert B_{r}(x^{0})\backslash E\rvert\\
 & \le Cr^{n}-C(\delta r)^{n}.
\end{align*}
Therefore the Lebesgue upper density \[ \limsup_{r\to0^{+}}\frac{\lvert E\cap B_{r}(x^{0})\rvert}{\lvert B_{r}(x^{0})\rvert}\le\limsup_{r\to0^{+}}\frac{r^{n}(1-\delta^{n})}{r^{n}}=1-\delta^{n}.
\] The Lebesgue upper density however is equal to 1 for a.e.~point in a set, which demonstrates that $\lvert E\cap\mathcal{D}\rvert=0$.
\end{proof}
\begin{thm}
[Porosity of the free boundary {\cite[Lemma 3.10]{MR2962060}}]\label{thm:porosity} Let $f\in L^{\infty}(\mathcal{D})$, $f\ge c>0$, and let $u$ be a solution to the obstacle problem in $\mathcal{D}$.
Then $\Gamma(u)$ is locally porous, hence $\lvert\Gamma(u)\rvert=0$.
\end{thm}

\begin{proof}
Let $K\Subset\mathcal{D}$, $x^{0}\in\Gamma(u)$ and $B_{2r}(x^{0})\subset K$.
By nondegeneracy, Lemma~\ref{lem:nondegeneracy}, and continuity of $u$, there exists $y\in\overline{B_{r/2}(x^{0})}$ such that \[ u(y)\ge\frac{c}{2n}\Big(\frac{r}{2}\Big)^{2}=\frac{c}{8n}r^{2}.
\] In particular $y\in\Omega(u)$, so $d\df\dist(y,\Gamma(u))$ is positive, and $d\le\lvert y-x^{0}\rvert\le r/2$.
Let $z\in\Gamma(u)$ be a point with $\lvert y-z\rvert=d$.
Then $\lvert z-x^{0}\rvert\le d+r/2$, so $B_{2d}(z)\subseteq B_{3d+r/2}(x^{0})\subseteq B_{2r}(x^{0})\subset\mathcal{D}$, and the sub-quadratic growth, Lemma~\ref{lem:sub-quadratic-growth}, together with continuity of $u$ gives \[ \frac{c}{8n}r^{2}\le u(y)\le\sup_{\overline{B_{d}(z)}}u\le C_{0}\lVert f\rVert_{L^{\infty}(\mathcal{D})}d^{2}, \] where $C_{0}=C_{0}(n)$ is the constant in Lemma~\ref{lem:sub-quadratic-growth}.
Hence $d\ge\delta r$ for \[ \delta\df\min\bigg(\sqrt{\frac{c}{8nC_{0}\lVert f\rVert_{L^{\infty}(\mathcal{D})}}},\frac{1}{2}\bigg), \] which only depends on $n$, $c$, and $\lVert f\rVert_{L^{\infty}(\mathcal{D})}$.
Consequently $B_{\delta r}(y)\cap\Gamma(u)=\varnothing$, and since $\lvert y-x^{0}\rvert\le r/2$ and $\delta r\le r/2$, \[ B_{\delta r}(y)\subseteq B_{r}(x^{0})\backslash\Gamma(u).
\]
\end{proof}

\begin{rem}
The results of this section so far, i.e.,~nondegeneracy, nondegeneracy of the gradient, and porosity of the free boundary, only require $f\in L^{\infty}(\mathcal{D})$ with $f\ge c>0$.
In particular, they rely on the sub-quadratic growth in Lemma~\ref{lem:sub-quadratic-growth} and the $W_{\text{loc}}^{2,p}$-regularity of $u$, but not on the optimal $C^{1,1}$-regularity, and hence no $C^{1,1}$-potential of $f$ is needed.
This changes in the next subsection, where uniform $C^{1,1}$-bounds of the rescaled solutions are essential.
\end{rem}

\subsection{Blow-Ups and Classification of Free Boundary Points}

Since we know that a solution $u$ grows quadratically away from the free boundary, we should be able to retain some information locally around a free boundary point $x^{0}\in\Gamma(u)$ by considering the rescaled functions
\begin{equation}
u_{x^{0},r}(x)\df\frac{u(x^{0}+rx)}{r^{2}}\label{eq:rescalings}
\end{equation}
and let $r>0$ tend to zero.
Throughout the rest of this section we assume, in addition to $f\ge c>0$, that $f$ is continuous and has a $C^{1,1}$-potential, so that $u\in C_{\text{loc}}^{1,1}(\mathcal{D})$ by Theorem~\ref{thm:c11-regularity}.
Since $\lvert D^{2}u_{x^{0},r}(x)\rvert =\lvert D^{2}u(x^{0}+rx)\rvert\le M<\infty$ for a.e.~$x$ in $B_{R/r}$ (for $0<R<\dist(x^{0},\partial\mathcal{D})$) due to the $C^{1,1}$-regularity of $u$, the Arzel\`{a}--Ascoli theorem implies that $\{u_{x^{0},r}\}_{r>0}$ is precompact in $C_{\text{loc}}^{1,\alpha}(\mathbb{R}^{n})$ for $\alpha\in[0,1)$.
Therefore there exists a subsequence $\{u_{x^{0},r_{j}}\}_{j=1}^{\infty}$, $r_{j}\to0^{+}$ as $j\to\infty$, and a function $u_{0}\in C_{\text{loc}}^{1,1}(\mathbb{R}^{n})$ such that $u_{x^{0},r_{j}}\to u  _{0}$ in $C_{\text{loc}}^{1,\alpha}(\mathbb{R}^{n})$.
Such a function, $u_{0}$, is called a \emph{blow-up limit} (or just a blow-up) of $u$ around $x^{0}$.
By the continuity of $f$ we see that, in the noncoincidence set,
\begin{equation}
\Delta u_{x^{0},r}(x)=f(x^{0}+rx)\to f(x^{0}),\qquad r\to0^{+}.\label{eq:blowup-value}
\end{equation}

In connection to these blow-ups, it turns out that the following two definitions are fruitful.
\begin{defn}
[Regular points] Let $u$ be a solution to the obstacle problem.
A point $x^{0}\in\Gamma(u)$ is called \emph{regular} if there exists a sequence $\{u_{x^{0},r_{j}}\}_{j=1}^{\infty}$ of rescalings of the form \eqref{eq:rescalings} that converges in $C_{\text{loc}}^{1,\alpha}(\mathbb{R}^{n})$ to a blow-up limit of the form \[ u_{0}(x)=\frac{f(x^{0})}{2}[(\mathbf{e}\cdot x)^{+}]^{2}, \] where $\mathbf{e}$ is a unit vector in $\mathbb{R}^{n}$.
\end{defn}

\begin{figure}[h]
\centering
\begin{tikzpicture}
\node[anchor=south east,inner sep=0] (image) at (0,0) {\includegraphics[width=0.8\textwidth]{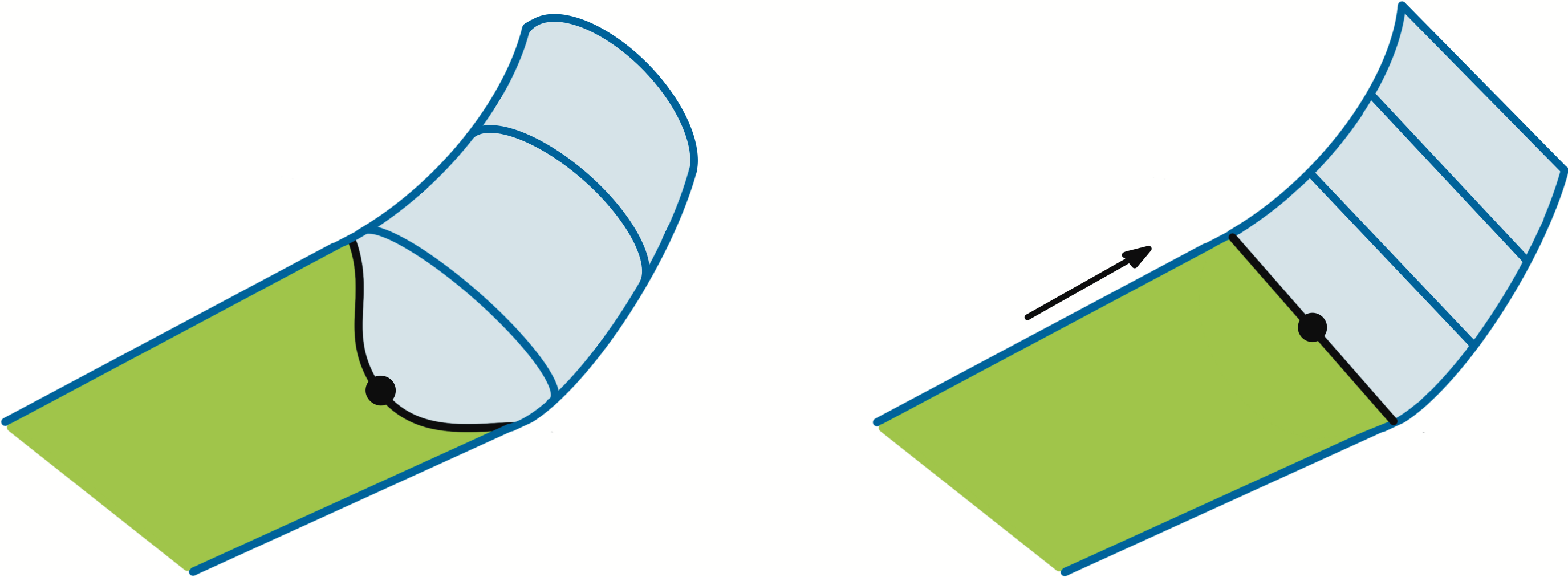}};
\begin{scope}[x={(image.south west)},y={(image.north east)}]
\node at (0.8,0.32) {$x^0$};
\node at (0.576,0.45) {$u$};
\node at (0.19,0.4) {$0$};
\node at (0.013,0.42) {$u_0$};
\node at (0.32,0.555) {$\mathbf{e}$};
\end{scope}
\node[below left] (copyright-text) {\tiny{\copyright Andreas Minne}};
\end{tikzpicture}\caption{Blow-up limit $u_{0}(x)=\frac{f(x^{0})}{2}[(\mathbf{e}\cdot x)^{+}]^{2}$
at a regular point $x^{0}$. }
\end{figure}
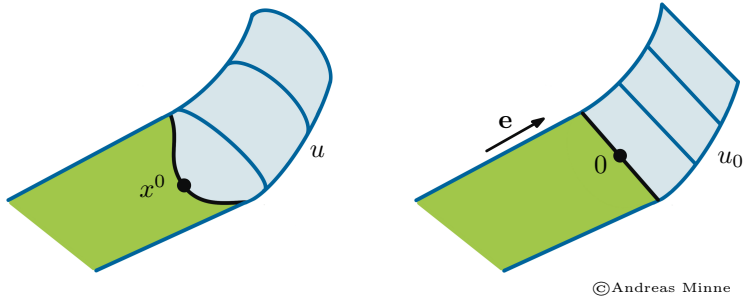
\begin{defn}
[Singular points] Let $u$ be a solution to the obstacle problem.
A point $x^{0}\in\Gamma(u)$ is called \emph{singular} if there exists a sequence $\{u_{x^{0},r_{j}}\}_{j=1}^{\infty}$ of rescalings of the form \eqref{eq:rescalings} that converges in $C_{\text{loc}}^{1,\alpha}(\mathbb{R}^{n})$ to a blow-up limit of the form \[ u_{0}(x)=\frac{f(x^{0})}{2}(Ax\cdot x), \] where $A$ is a nonnegative definite matrix such that $\trace A=1$.
\end{defn}

\begin{figure}[h]
\centering
\begin{tikzpicture}
\node[anchor=south east,inner sep=0] (image) at (0,0) {\includegraphics[width=0.8\textwidth]{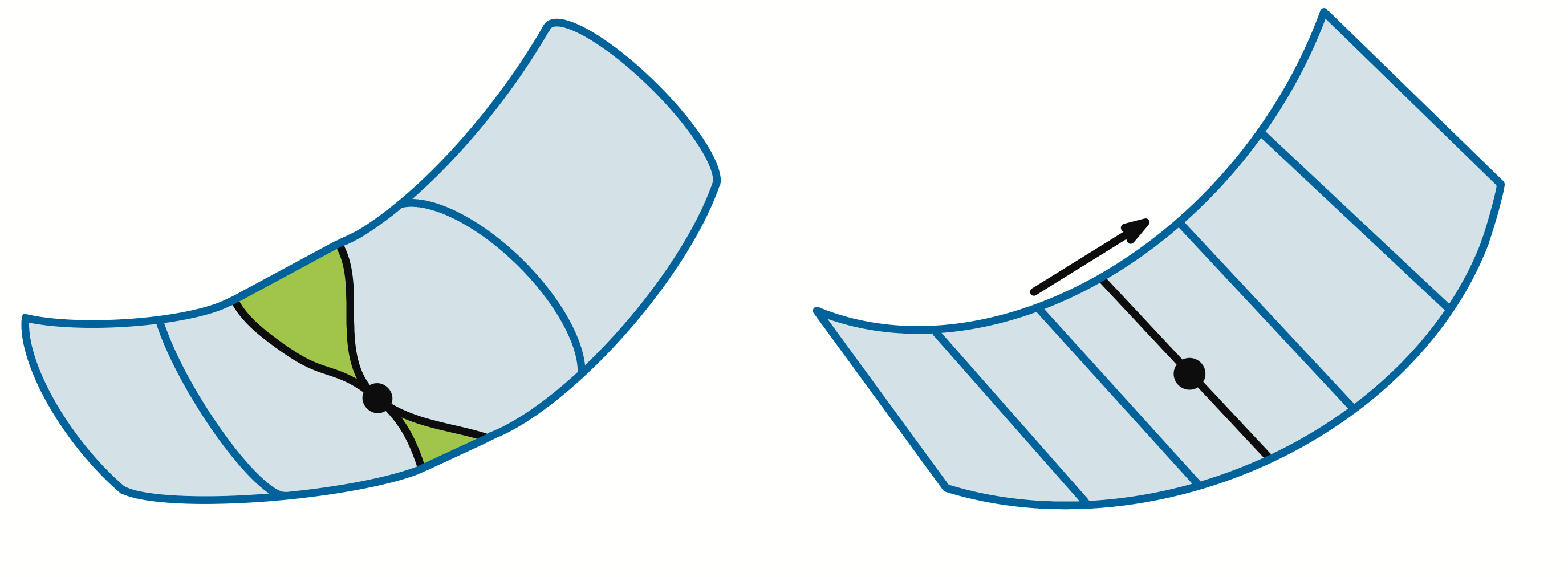}};
\begin{scope}[x={(image.south west)},y={(image.north east)}]
\node at (0.79,0.28) {$x^0$};
\node at (0.56,0.45) {$u$};
\node at (0.22,0.38) {$0$};
\node at (0.04,0.42) {$u_0$};
\node at (0.32,0.60) {$\mathbf{e}$};
\end{scope}
\node[below left] (copyright-text) {\tiny{\copyright Andreas Minne}};
\end{tikzpicture}\caption{Blow-up limit $u_{0}(x)=\frac{f(x^{0})}{2}(Ax\cdot x)$ at a singular
point $x^{0}$. }
\end{figure}
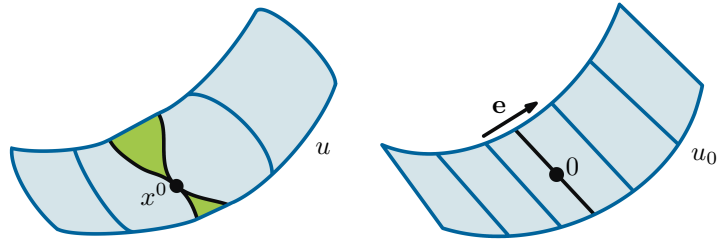

\begin{rem}
Note that the Laplacian of the blow-up limits above equal $f(x^{0})$ as alluded to in \eqref{eq:blowup-value}.
\end{rem}

\begin{example}
Let $f\equiv 1$.
Each of the following is a global solution to the obstacle problem on $\mathbb{R}^{n}$:
	\begin{enumerate}[label=(\roman{enumi})]
		\item
$u(x) = \frac{1}{2}(x_{1}^{+})^{2}$, where $\Gamma(u) = \{x_{1}=0\}$ consists entirely of regular points;
		\item
$u(x) = \frac{\lvert x\rvert^{2}}{2n}$, where $\Gamma(u) = \{0\}$ is a single singular point;
		\item
$u(x) = \frac{1}{2}x_{1}^{2}$, where $\Gamma(u) = \{x_{1}=0\}$ consists entirely of singular points.
	\end{enumerate}
Note that all blow-up limits coincide with the solutions themselves, since the solutions are homogeneous of degree two.
For (i), $\mathbf{e}$ points in the $x_1$-direction; for (ii), $A = I/n$; and for (iii), $A_{11}=1$ while all the other entries of $A$ are zero.
\end{example}

At the moment it is unclear that the definitions above are useful since different subsequences might lead to different blow-ups.
In fact, the following theorem by Blank \cite{Bla01} shows that we have to be careful.
\begin{thm}
[{\cite[Theorem 0.2]{Bla01}}] Let $\sigma$ be any modulus of continuity which is not Dini continuous, i.e.,
\[
\int_{0}^{\epsilon}\frac{\sigma(r)}{r}\,dr=\infty
\] 
for any $\epsilon>0$.
We can find an obstacle $\psi$ such that $\Delta\psi$ has $\sigma$ as its modulus of continuity, and such that the set between the obstacle and membrane spins an infinite number of times around a free boundary point where this set has density $1/2$.
\end{thm}

\begin{rem}
Note that the class of Dini continuous functions is strictly contained in the class of functions with $C^{1,1}$-potential.
\end{rem}

At this point it might be good to step back for a moment and reflect on the differences between the analysis of the solution $u$ and that of the free boundary.
We know that to obtain optimal regularity, it is enough that $f$ has a $C^{1,1}$-potential, while just to have uniqueness of blow-ups around a free boundary point, Dini continuity is needed.
The moral of the story is that the regularity of solutions and that of the free boundary should be treated separately, as the sharp assumptions in the two settings differ.

Before we proceed to the last theorem --- namely that there are no other types of blow-ups, i.e.,~$\Gamma(u)$ consists only of regular and singular points, and that the \emph{type} of blow-up is unique --- we will prove some basic properties of blow-up sequences.
\begin{lem}
[{\cite[Proposition 3.17]{MR2962060}}]\label{lem:blow-up-properties} Let $u$ be a solution to the obstacle problem in $\mathcal{D}$ with a Dini continuous function $f\ge c > 0$, a point $x^{0}\in\Gamma(u)$, and $\{u_{j}\}_{j=1}^{\infty}\df\{u_{x^{0},r_{j}}\}_{j=1}^{\infty}$ a sequence of rescalings of the form \eqref{eq:rescalings} that converges in $C_{\text{loc}}^{1,\alpha}(\mathbb{R}^{n})$ to a blow-up limit $u_{0}\in C_{\text{loc}}^{1,1}(\mathbb{R}^{n})$.
Then:
\begin{enumerate}[label=(\roman{enumi})]
\item
$u_{0}(y^{0})>0\implies u_{j}>0$ in $B_{\delta}(y^{0}),$ for $j$ large enough.
\item
$\lvert\nabla u_{0}(y^{0})\rvert >0\implies\lvert\nabla u_{j}\rvert >0$ in $B_{\delta}(y^0),$ for $j$ large enough.
\item
$\nabla u_{0}=0$ in $B_{\delta}(y^{0})\implies\nabla u_{j}=0$ in $B_{\delta/2}(y^0)$ for $j$ large enough.
\item
$u_{0}$ solves the obstacle problem in $\mathbb{R}^{n}$.
\item
$u_{j}\to u_{0}$ strongly in $W_{\text{loc}}^{2,p}(\mathbb{R}^{n})$ for any $p\in(1,\infty)$.
\end{enumerate}
\end{lem}

\begin{proof}
\textbf{(i) and (ii)}: These immediately follow from the $C_{\text{loc}}^{1,\alpha}$-convergence.

\noindent \textbf{(iii)}: Suppose that there exists a subsequence $\{u_{j_{k}}\}_{k=1}^{\infty}$ such that $\lvert\nabla u_{j_{k}}(y^{k})\rvert >0$ for some $y^{k}\in B_{\delta/2}(y^{0})$ but $\lvert\nabla u_{0}\rvert =0$ in $B_{\delta}(y^{0})$. 
The points $y^k$ lie in $\overline{\Omega(u_{j_k})}$, and since the rescaled functions $u_{j_{k}}$ also solve the obstacle problem in balls around the origin of radius $R/r_{j_{k}}$ for $0<R<\dist(x^{0},\partial\mathcal{D})$, by nondegeneracy of the gradient (Corollary~\ref{cor:nondegeneracy-gradient}), \[ \sup_{B_{3\delta/4}(y^{0})}\lvert\nabla u_{j_{k}}\rvert\ge\sup_{B_{\delta/4}(y^{k})}\lvert\nabla u_{j_{k}}\rvert\ge\frac{c\delta}{8n}.
\] By letting $k\to\infty$, we get $\sup_{B_{3\delta/4}(y^{0})}\lvert\nabla u_{0}\rvert\ge c\delta/(8n)$, a contradiction.

\noindent \textbf{(iv)}: The blow-up $u_{0}$ is in $W_{\text{loc}}^{2,\infty}(\mathbb{R}^{n})$.
Therefore it is enough to check that the equation \[
\begin{cases}
\Delta u_{0}=f(x^{0})\chi_{\{u_{0}>0\}},\\
\phantom{\Delta}u_{0}\ge0,
\end{cases}
\] is satisfied almost everywhere in $\mathbb{R}^{n}$: Since $D^{2}u_{0}=0$ a.e.~in $\Lambda(u_{0})$, also $\Delta u_{0}=\trace D^{2}u_{0}=0$ a.e.~in $\Lambda(u_{0})$.
For $y^{0}$ in $\Omega(u_{0})$, by (i) there is a $\delta>0$ and $j_{0}\in\mathbb{N}$ such that for $j\ge j_{0}$, $B_{\delta}(y^{0})\subseteq\Omega(u_{j})$, hence $\Delta u_{j}(x)=f(x^{0}+r_{j}x)$ in $B_{\delta}(y^{0})$.
Also, since $\lvert f(x^{0}+r_{j}x)-f(x^{0})\rvert\le\sigma(r_{j}\lvert x\rvert )\le\sigma\big(r_{j}(\lvert y^{0}\rvert+\delta)\big)\to0$ as $j\to\infty$ for $\sigma$ a modulus of continuity of $f$, $\Delta u_{j}(x)$ converges uniformly to $f(x^{0})$ in $B_{\delta}(y^{0})$, and classical results then imply that $\Delta u_{0}$ exists and is equal to $f(x^{0})$ in this set.
In summary, \[ \Delta u_{0}=f(x^{0})\chi_{\{u_{0}>0\}}.
\]

\noindent \textbf{(v)}: Fix a compact set $K\subset\mathbb{R}^{n}$ and $p\in(1,\infty)$.
By the $C^{1,1}$-regularity, $\lvert D^{2}u_{j}\rvert\le M$ a.e.~in $K$ for $j$ large, and also $\lvert D^{2}u_{0}\rvert\le M$.
We split $\mathbb{R}^{n}$ into the three disjoint sets $\Omega(u_{0})$, $\interior(\Lambda(u_{0}))$, and $\Gamma(u_{0})$, and treat them separately.

\emph{The free boundary $\Gamma(u_{0})$.} By (iv), $u_{0}$ solves the obstacle problem with the constant right-hand side $f(x^{0})\ge c>0$, so $\lvert\Gamma(u_{0})\rvert=0$ by Theorem~\ref{thm:porosity}.
Hence, given $\epsilon>0$, there are compact sets $K_{1}\subset\Omega(u_{0})$ and $K_{2}\subset\interior(\Lambda(u_{0}))$ such that $\lvert K\backslash(K_{1}\cup K_{2})\rvert<\epsilon$.

\emph{The noncoincidence set $\Omega(u_{0})$.} For $y^{0}\in K_{1}$, by continuity of $u_{0}$ and by (i), there is $\delta>0$ such that $u_{0}>0$ and $u_{j}>0$ in $B_{\delta}(y^{0})$ for $j$ large.
The difference $w_{j}\df u_{j}-u_{0}$ then satisfies, in $B_{\delta}(y^{0})$,
\[
\lvert\Delta w_{j}(x)\rvert=\lvert f(x^{0}+r_{j}x)-f(x^{0})\rvert\le\sigma\big(r_{j}(\lvert y^{0}\rvert+\delta)\big)\to0,
\]
where $\sigma$ is a modulus of continuity of $f$, and also $w_{j}\to0$ uniformly in $B_{\delta}(y^{0})$.
The Calder\'{o}n--Zygmund estimates in Theorem~\ref{thm:cz-estimates} applied to $w_{j}$ give
\[
\lVert D^{2}w_{j}\rVert_{L^{p}(B_{\delta/2}(y^{0}))}\le C\left(\lVert w_{j}\rVert_{L^{1}(B_{\delta}(y^{0}))}+\lVert\Delta w_{j}\rVert_{L^{p}(B_{\delta}(y^{0}))}\right)\to0.
\]
By compactness, $K_{1}$ is covered by finitely many such balls $B_{\delta/2}(y^{0})$, so $D^{2}u_{j}\to D^{2}u_{0}$ in $L^{p}(K_{1})$.

\emph{The interior of the coincidence set.} For $y^{0}\in K_{2}$, choose $\delta>0$ such that $B_{\delta}(y^{0})\subseteq\interior(\Lambda(u_{0}))$.
Then $\nabla u_{0}=0$ in $B_{\delta}(y^{0})$, so by (iii), $\nabla u_{j}=0$ in $B_{\delta/2}(y^{0})$ for $j$ large, and hence $D^{2}u_{j}=0=D^{2}u_{0}$ there.
Again by compactness, $D^{2}u_{j}=D^{2}u_{0}$ a.e.~in $K_{2}$ for $j$ large.

\emph{Conclusion.} Combining the three parts,
\begin{align*}
\limsup_{j\to\infty}\int_{K}\lvert D^{2}u_{j}-D^{2}u_{0}\rvert^{p}\,dx & \le\limsup_{j\to\infty}\int_{K_{1}}\lvert D^{2}u_{j}-D^{2}u_{0}\rvert^{p}\,dx+(2M)^{p}\epsilon\\
 & =(2M)^{p}\epsilon.
\end{align*}
Since $\epsilon>0$ is arbitrary, $D^{2}u_{j}\to D^{2}u_{0}$ in $L^{p}(K)$, and together with the $C^{1,\alpha}_{\text{loc}}$-convergence this shows that $u_{j}\to u_{0}$ in $W_{\text{loc}}^{2,p}(\mathbb{R}^{n})$.
\end{proof}
The next theorem is one of the most famous tools within the field of FBPs: the Alt--Caffarelli--Friedman monotonicity formula \cite{MR732100}.
The proof of this theorem can be found in \cite{MR2962060}.
\begin{thm}
[{\cite[Lemma 5.1]{MR732100}}]\label{thm:ACF} Let $u_{1}$, $u_{2}\in C(B_{1})\cap W^{1,2}(B_{1})$ be two functions such that, for $i=1,2$, \[ u_{i}\ge0,\qquad\Delta u_{i}\ge0,\qquad u_{1}\cdot u_{2}=0,\qquad\text{in }B_{1}.
\] Then \[ \Phi(r,u_{1},u_{2})\df\frac{1}{r^{4}}\int_{B_{r}}\frac{\lvert\nabla u_{1}\rvert^{2}}{\lvert x\rvert^{n-2}}\,dx\int_{B_{r}}\frac{\lvert\nabla u_{2}\rvert^{2}}{\lvert x\rvert^{n-2}}\,dx \] is nondecreasing for $0<r<1$.
If $\Phi(r_{1})=\Phi(r_{2})$ for some $r_1$ and $r_2$ such that $0<r_{1}<r_{2}<1$, then:
\begin{enumerate}[label=(\roman{enumi})]
\item
at least one of $u_{1}$ or $u_{2}$ is zero in $B_{r_{2}}$, or
\item
there exist constants $k_{i}>0$, $i=1$, $2$ and a unit vector $\mathbf{e}$ such that \[ u_{1}=k_{1}(\mathbf{e}\cdot x)^{+},\qquad u_{2}=k_{2}(\mathbf{e}\cdot x)^{-},\qquad\text{in }B_{r_{2}}.
\]
\end{enumerate}
\end{thm}

\begin{rem}
The continuity assumption can in fact be dropped, see for instance \cite{vel14}.
\end{rem}

To get some familiarity with how such a monotonicity formula can be used, we will use it to prove the final theorem of this text, which shows that all the free boundary points are either regular or singular, and that they are mutually exclusive, i.e.,~a free boundary point cannot be both regular and singular.
For simplicity we assume that $f$ is identically equal to one.
\begin{thm}[Classification of free boundary points]\label{thm:classification}
Let $u$ be a solution to the obstacle problem in $\mathcal{D}$ with $f(x)\equiv 1$, a point $x^{0}\in\Gamma(u)$, and $\{u_{j}\}_{j=1}^{\infty}\df\{u_{x^{0},r_{j}}\}_{j=1}^{\infty}$ a sequence of rescalings of the form \eqref{eq:rescalings} that converges in $C_{\text{loc}}^{1,\alpha}(\mathbb{R}^{n})$ to a blow-up limit $u_{0}$.
Then either \[ u_{0}(x)=\frac{1}{2}[(\mathbf{e}\cdot x)^{+}]^{2},\qquad\text{or}\qquad u_{0}(x)=\frac{1}{2}(Ax\cdot x), \] where $\mathbf{e}$ is a unit vector, and $A$ is a nonnegative definite matrix such that $\trace A=1$.
Moreover, the type of blow-up is unique.
\end{thm}

\begin{proof}
Let $\omega\in\partial B_{1}$ be a unit vector, $r>0$ so that $B_{r}(x^{0})\Subset\mathcal{D}$ and \[ \phi_{\omega}(r,u,x^{0})\df\frac{1}{r^{4}}\int_{B_{r}(x^{0})}\frac{\lvert\nabla(\partial_{\omega}u)^{+}\rvert^{2}}{\lvert x-x^{0}\rvert^{n-2}}\,dx\int_{B_{r}(x^{0})}\frac{\lvert\nabla(\partial_{\omega}u)^{-}\rvert^{2}}{\lvert x-x^{0}\rvert^{n-2}}\,dx.
\] Since $u_{x^{0},r_{j}}\to u_{0}$ strongly in $W_{\text{loc}}^{2,p}(\mathbb{R}^{n})$ by Lemma~\ref{lem:blow-up-properties}(v), $\phi_{\omega}(r,u_{0},0)=\lim_{j\to\infty}\phi_{\omega}(r,u_{x^{0},r_{j}},0)$.
We also have the rescaling property $\phi_{\omega}(r,u_{x^{0},r_{j}},0)=\phi_{\omega}(rr_{j},u,x^{0})$.
Let us now check the conditions for the Alt--Caffarelli--Friedman formula: By definition, $(\partial_{\omega}u)^{+}\cdot(\partial_{\omega}u)^{-}=0$ , and $(\partial_{\omega}u)^{+}\ge0$, $(\partial_{\omega}u)^{-}\ge0$ in $\mathcal{D}$.
Since $\nabla u=0$ on $\Lambda(u)$, we have $\{y\in\mathcal{D}:\partial_{\omega}u(y)\ne0\}\subseteq\Omega(u)$.
Hence, as $\partial_{\omega}u$ is continuous, for $x\in\{y\in\mathcal{D}:\partial_{\omega}u(y)>0\}$ there exists a ball $B_{\rho}(x)\subseteq\{y\in\mathcal{D}:\partial_{\omega}u(y)>0\}\subseteq\Omega(u)$, $\rho>0$.
Therefore $\Delta u=1$ in $B_{\rho}(x)$, hence $\Delta\partial_{\omega}u=\partial_{\omega}\Delta u=0$, and similar for $x\in\{y\in\mathcal{D}:\partial_{\omega}u(y)<0\}$.
For $x\in\{y\in\mathcal{D}:\partial_{\omega}u(y)=0\}$, we have \[ (\partial_{\omega}u)^{\pm}(x)=0\le\fint_{B_{\rho}(x)}(\partial_{\omega}u)^{\pm}(y)\,dy, \] where $\fint_{B}\df\frac{1}{\lvert B\rvert }\int$ is the mean value integral.
In conclusion, both $(\partial_{\omega}u)^{+}$ and $(\partial_{\omega}u)^{-}$ are subharmonic.
This means that we can apply the monotonicity formula.
The limit $\lim_{r\to0^{+}}\phi_{\omega}(r,u,x^{0})\fd L_{\omega}$ exists since $\phi_{\omega}(r,u,x^{0})$ is nondecreasing in $r$ and bounded from below by zero.
It is also finite: since $\lvert\nabla(\partial_{\omega}u)^{\pm}\rvert\le\lvert D^{2}u\rvert$ a.e.~and $\int_{B_{r}(x^{0})}\lvert x-x^{0}\rvert^{2-n}\,dx=C(n)r^{2}$, the $C^{1,1}$-regularity of $u$ gives \[ \phi_{\omega}(r,u,x^{0})\le\frac{1}{r^{4}}\Big(C(n)r^{2}\lVert D^{2}u\rVert_{L^{\infty}(B_{r}(x^{0}))}^{2}\Big)^{2}=C(n)^{2}\lVert D^{2}u\rVert_{L^{\infty}(B_{r}(x^{0}))}^{4}.
\] Moreover, \[ \phi_{\omega}(r,u_{0},0)=\lim_{j\to\infty}\phi_{\omega}(r,u_{x^{0},r_{j}},0)=\lim_{j\to\infty}\phi_{\omega}(rr_{j},u,x^{0})=L_{\omega}\ge0.
\] Therefore $\phi_{\omega}(r,u_{0},0)$ is constant and the case of equality in the Alt--Caffarelli--Friedman monotonicity formula gives exactly one of the following possibilities:
\begin{enumerate}[label=(\roman{enumi})]
\item
Either $(\partial_{\omega}u_{0})^{+}\equiv0$ or $(\partial_{\omega}u_{0})^{-}\equiv0$ in $\mathbb{R}^{n}$,
\item
$(\partial_{\omega}u_{0})^{+}\equiv k_{1}(\mathbf{e}\cdot x)^{+}$ and $(\partial_{\omega}u_{0})^{-}\equiv k_{2}(\mathbf{e}\cdot x)^{-}$.
\end{enumerate}
Let us consider two cases:

\noindent\textbf{Case (a)} ($\lvert\{x\in\mathbb{R}^{n}:u_{0}(x)=0\}\rvert >0$).
If the measure is nonzero, (i) holds for any direction $\omega$.
Indeed, in (ii) the function $\partial_{\omega}u_{0}=k_{1}(\mathbf{e}\cdot x)^{+}-k_{2}(\mathbf{e}\cdot x)^{-}$ vanishes only on the hyperplane $\{\mathbf{e}\cdot x=0\}$, whereas $\nabla u_{0}=0$ on the set $\{u_{0}=0\}$ of positive measure.
This shows that $\partial_{\omega}u_0$ has a sign for every $\omega$; by Lemma~\ref{lem:direction-sign} below, $u_{0}$ is a monotone function of one variable, say $u_{0}(x)=v(x_{1})$ without loss of generality, where $v$ is nondecreasing in $x_{1}$.

Now we claim that $0\in\Gamma(u_{0})$.
If $0\in\Omega(u_{0})$, then $u_{x^{0},r_{j}}>0$ in $B_{\delta}$ for some $\delta>0$ and $j$ large by Lemma~\ref{lem:blow-up-properties}(i), contradicting that $x^{0}\in\Gamma(u)$.
Moreover, if $0\in\interior(\Lambda(u_{0}))$, then $\nabla u_{x^{0},r_{j}}\equiv0$, in $B_{\delta}$ for some $\delta>0$ by Lemma~\ref{lem:blow-up-properties}(iii), again contradicting the fact that $x^{0}\in\Gamma(u)$.
Therefore $0\in\Gamma(u_{0})$.

Since $v$ is nondecreasing, $v\equiv0$ in the halfspace $\{x\in\mathbb{R}^{n},x_{1}\le0\}$ due to being above the (zero) obstacle, and $\Delta u_{0}=\partial_{11}v=1$ in $\{x\in\mathbb{R}^{n},x_{1}>0\}$ since $0\in \Gamma(u_0)$.
This implies that $v=x_{1}^{2}/2+bx_1+c$ for $x_{1}>0$, and since $v(0)=v'(0)=0$ we get $b=c=0$.
Hence $v=x_{1}^{2}/2$ in the positive halfspace, so $u_{0}=(x_{1}^{+})^{2}/2$ in $\mathbb{R}^{n}$.

\noindent\textbf{Case (b)} ($\lvert\{x\in\mathbb{R}^{n}:u_{0}(x)=0\}\rvert =0$).
In this case, $\Delta u_{0}=1$ almost everywhere in $\mathbb{R}^{n}$, and so from the classical regularity theory we get that $u_{0}$ is analytic, and that \[ \Delta\partial_{ij}u_{0}=\partial_{ij}\Delta u_{0}=0.
\] Also, for some small $\delta>0$, $\lvert\partial_{ij}u_{0}\rvert\le\lVert D^{2}u\rVert_{L^{\infty}(B_{\delta}(x^{0}))}<\infty$.
Therefore the Liouville theorem gives $D^{2}u_{0}=A$ for some symmetric matrix $A$.
From this, \[ u_{0}(x)=\frac{1}{2}(Ax\cdot x)+b\cdot x+c,\qquad b\in\mathbb{R}^{n},c\in\mathbb{R}.
\] The conditions $u_{0}(0)=\lvert\nabla u_{0}(0)\rvert =0$ again show that $\lvert b\rvert=c=0$.
Moreover, $\Delta u_{0}=\trace A=1$, and $u_{0}\ge0$ shows that $A$ has to be nonnegative definite.

For uniqueness of type of blow-up limit: Assume that there exists a blow-up $u_{0}$ of the form $u_{0}=\frac{1}{2}[(\mathbf{e}\cdot x)^{+}]^{2}$.
Then \[ 0=\phi_{\omega}(r,u_{0},0)=\lim_{j\to\infty}\phi_{\omega}(rr_{j},u,x^{0})=L_{\omega} \] for all unit vectors $\omega$.
If $v_{0}$ is another blow-up, obtained along a sequence $s_{j}\to0^{+}$, then \[ \phi_{\omega}(r,v_{0},0)=\lim_{j\to\infty}\phi_{\omega}(rs_{j},u,x^{0})=L_{\omega}=0 \] from which it can be deduced that $v_{0}$ cannot be of the form $\frac{1}{2}(Ax\cdot x)$ because otherwise there would be a direction $\omega$ such that $\phi_{\omega}(r,v_{0},0)>0$, and vice versa if we started out with a blow-up limit $u_{0}$ of the form $\frac{1}{2}(Ax\cdot x)$.
Hence the type of blow-up is unique.
\end{proof}
\begin{lem}
\label{lem:direction-sign} If $u\in C^{1}(\mathbb{R}^{n})$ and $\partial_{\omega}u$
has a sign for every direction $\omega$, then $u$ is a monotone function of one variable.
\end{lem}

\begin{proof}
If $\nabla u(x)=\lvert\nabla u(x)\rvert\mathbf{e}$ for some fixed direction $\mathbf{e}$, then $\partial_{\mathbf{e}'}u=\lvert\nabla u\rvert\mathbf{e}\cdot\mathbf{e}'=0$ for every $\mathbf{e}'\perp\mathbf{e}$ and all $x\in\mathbb{R}^{n}$, which implies that $u(x)=v(\mathbf{e}\cdot x)$, and $v$ is nondecreasing.

If there are two different points $x^{1}$, $x^{2}$ and directions $\mathbf{e}$, $\mathbf{e}'$, $\mathbf{e}\ne\mathbf{e}'$ such that $0\ne\lvert\nabla u(x^{1})\rvert\mathbf{e}$ and $0\ne\lvert\nabla u(x^{2})\rvert\mathbf{e}'$, then there exists a unit vector $\omega$ for which $\mathbf{e}\cdot\omega>0$ and $\mathbf{e}'\cdot\omega<0$.
Consequently $\partial_{\omega}u(x^{1})=\lvert\nabla u(x^{1})\rvert\mathbf{e}\cdot\omega>0$ and $\partial_{\omega}u(x^{2})=\lvert\nabla u(x^{2})\rvert\mathbf{e}'\cdot\omega<0$ which would contradict that $\partial_{\omega}u$ has a sign.
\end{proof}

\subsection*{Further Reading}

The blow-up classification proved in Theorem~\ref{thm:classification} is the first step toward the classical regularity theorem of Caffarelli~\cite{MR0454350}: in a neighborhood of a regular free boundary point, $\Gamma(u)$ is a $C^{1,\alpha}$ hypersurface.
The interested reader is referred to the monograph~\cite{MR2962060} for this and much more.

\subsubsection*{Acknowledgements}

It is my pleasure to express my appreciation to Professor Edgard Pimentel for inviting me to the Rio Summer School on PDEs 2020, and for suggesting that I write this text for Coimbra Mathematical Texts.
I am also very grateful to the anonymous referee, as well as Anders Furufors, Malte Litsg\r{a}rd and David Tewodrose, for taking the time to read this text and for providing valuable feedback and corrections.
Finally, I acknowledge the use of AI tools (Claude and ChatGPT) for editing and proofreading.
I remain responsible for the final content.

\bibliographystyle{amsalpha}
\bibliography{References}

\end{document}